%% file: CLInstGen_TR.tex
\documentclass[3p,authoryear]{elsarticle2}

\usepackage{etoolbox}
\newbool{TR}
\setbool{TR}{true}

\AtBeginEnvironment{table}{}

\usepackage[T1]{fontenc}
\usepackage{pifont}
\usepackage[english]{babel}
\usepackage{marvosym}
\usepackage{natbib}
\usepackage{pdflscape}
\normalsize

\usepackage{fontawesome5} % nice Symbols
\usepackage{longtable}
\usepackage{pdflscape} % Allows for landscape pages
\usepackage{geometry} % Adjusts page size if needed
\usepackage{multirow} %enthält multirow-Befehl für Tabellen
\usepackage{booktabs} % Für bessere Tabellenlinien
\usepackage{rotating} % Ermöglicht das Drehen von Tabellen
\usepackage{geometry} % Passt die Seitengröße an, falls erforderlich
\usepackage{subcaption}
\usepackage{caption}
\usepackage{array}
\usepackage{mathrsfs,amssymb,amsthm}
\usepackage[intlimits]{empheq}
\DeclareMathAlphabet{\mathcal}{OMS}{cmsy}{m}{n}
\usepackage{mathtools} %enthält multilined-Befehl für Align-Umgebungen
\numberwithin{equation}{section}
\allowdisplaybreaks %erlaubt Seitenumbrüche auch in alignat-Umgebungen

\usepackage{xcolor}
\definecolor{rptured}{RGB}{227,27,76}
\definecolor{rptublgrey}{RGB}{80,114,137}

\usepackage{tikz}
\usetikzlibrary{arrows, shadings, shapes, decorations.pathmorphing,patterns, decorations.pathreplacing} % for triangle45 arrows, decorations.pathreplacing for curly bracket
\usepackage{comment}
\usepackage{todonotes}

\usepackage{afterpage}

\ifbool{TR}
{
	\input mytableTR
	\input myfigureTR
}
{
	\input mytableInforms
	\input myfigureInforms
}

\input{CLIG_Commands}

\begin{document}
	\begin{frontmatter}
		\title{CLInstGen: An Instance Generation Framework for Two-Tier City Logistics}
		% Benchmark Instance Generation for Two-Tier City Logistics
		\author[rptu]{Julia Lange\corref{cor1}}
		\cortext[cor1]{Corresponding author.}
		\ead{julia.lange@rptu.de}
		\author[ku]{Johannes Gückel}
		\ead{jgueckel@ku.de}
		\author[ku]{Pirmin Fontaine}
		\ead{pirmin.fontaine@ku.de}
		%\author[rptu]{Timo Gschwind\corref{cor2}}
		\author[rptu]{Timo Gschwind}
		\ead{gschwind@rptu.de}
		\address[rptu]{Chair of Logistics, School of Business and Economics, RPTU Kaiserslautern-Landau, 67663 Kaiserslautern, Germany.}
		\address[ku]{Chair of Logistics and Operations Analytics, Ingolstadt School of Management, Catholic University Eichstätt Ingolstadt, 85049 Ingolstadt, Germany}
		
		\ifbool{TR}
		{\journal{L-2026-02}}
		{\journal{}}
		
		\begin{abstract}
			City logistics addresses the design and management of efficient freight transportation in urban areas. Two-tier networks, where goods are moved by different vehicle fleets from distribution centers to customers using intermediate consolidation facilities, are particularly effective. Although studied from routing, location, and network design perspectives, the research progress is slowed by the lack of standardized, state-of-the-art benchmark data. Many studies rely on outdated and artificial or case-specific instances, limiting reproducibility, comparability, and generality of results. 
			To address this gap, we identify general problem characteristics for two-tier city logistics, provide an instance generation framework CLInstGen that is applicable to a wide range of problem settings, and propose the first benchmark library CLLIB. The library and the CLInstGen code are open-source, allowing researchers to easily use, adapt, and generate instances for their specific applications. CLLIB is based on real-world geographies of urban regions and provides a curated collection of standardized instances for direct benchmarking of models and methods as well as statistical evaluation. 
			A computational study with two different problem settings, namely tactical service network design on two tiers and vehicle routing with simultaneous pickups and deliveries on a single tier, confirms the versatility of both the generator and the library, and provides first numerical and managerial insights on the CLLIB benchmark instances.
		\end{abstract}

		
	\end{frontmatter}

%\newpage
\section{Introduction}
\label{sec:Intro}

\input{CLIG_1-Intro}

\section{Characterization of Problems in City Logistics}
\label{sec:InstDes}

\input{CLIG_3-InstDes}

\section{Instance Generation Framework for City Logistics}
\label{sec:InstGen}

\input{CLIG_4-InstGen}

\section{Computational Study on Applicability}
\label{sec:Study}

\input{CLIG_6-CompStud}

\section{Conclusion}
\label{sec:Concl}

\input{CLIG_7-Concl}

\ifbool{TR}
{\section*{Acknowledgments}

This research was supported by the German Research Foundation (project number: ROCOCO FO1468/1-1) and the Dr.~Werner Jackstädt-Stiftung. %
This support is gratefully acknowledged.
}
{}

\bibliography{CLIG_Paper.bib}

\ifbool{TR}
{\bibliographystyle{natbib}}
{\bibliographystyle{abbrvnat}}

\newpage
\appendix
\section{Maps of Urban Regions and Locations Used}
\label{app:Maps}
\input{CLIG_App_maps}

\section{CLLIB Sub-Libraries}
\label{app:SubLibr}
\input{CLIG_App_sublibraries}

\section{Details on Computational Studies}
\label{app:Model}
\input{CLIG_App_model}

\end{document}

%% file: mytableTR.tex
% Command mytable to use in our Technical Reports
%
% 1. Parameter (optional!): note under the table
% 2. Parameter: caption
% 3. Parameter: label
%	4. Parameter: Preference of position, e.g. hbt
% 5. Parameter: contents
%
%Standard use may be: \mytable[]{Caption of table}{tab:Label}{hbt}{\begin{tabu}...\{end{tabu}}

\newcommand{\mytable}[5][]{

\begin{table}[#4]
\caption{#2}
\label{#3} 
	
#5

#1
\end{table}

}

%% file: myfigureTR.tex
% Command myfigure to use in our Technical Reports
%
% 1. Parameter (optional!): note under the figure
% 2. Parameter: caption
% 3. Parameter: label
%	4. Parameter: Preference of position, e.g. hbt
% 5. Parameter: contents
%
%Standard use may be: \myfigure[]{Caption of figure}{fig:Label}{hbt}{\begin{tikzpicture}...\{end{tikzpicture}}

\newcommand{\myfigure}[5][]{
\begin{figure}[#4]
#5
\vspace{-0.25cm}
\caption{#2}
\label{#3}
#1
\end{figure}
}

%% file: mytableInforms.tex
% Command myfigure to use for Informs Journals
%
% 1. Parameter (optional!): note under the table
% 2. Parameter: caption
% 3. Parameter: label
% 4. Parameter: Preference of position, e.g. hbt
% 5. Parameter: contents
%
%Standard use may be: \mytable[]{Caption of table}{tab:Label}{hbt}{\begin{tabu}...\{end{tabu}}

\newcommand{\mytable}[5][]{
\begin{table}[#4]
\TABLE
{#2\label{#3}}
{#5}
{#1}
\end{table}
} 

%% file: myfigureInforms.tex
% Command myfigure to use for Informs Journals
%
% 1. Parameter (optional!): note under the figure
% 2. Parameter: caption
% 3. Parameter: label
%	4. Parameter: Preference of position, e.g. hbt
% 5. Parameter: contents
%
%Standard use may be: \myfigure[]{Caption of figure}{fig:Label}{hbt}{\begin{tikzpicture}...\{end{tikzpicture}}
%
%Attention: No \begin{center}...\end{center} in contents allowed!

\newcommand{\myfigure}[5][]{
\begin{figure}[#4]
\FIGURE
{#5}
{#2 \label{#3}}
{#1}
\end{figure}
}

%% file: CLIG_Commands.tex
\usepackage{glossaries}
\usepackage{pifont}
\usepackage{adjustbox}

\makeglossaries
\newacronym{dc}{DC}{Distribution Center}
\newacronym{lsp}{LSP}{Logistics Service Provider}
\newacronym{tw}{TW}{time window}
\newacronym{2elrp}{2E-LRP}{Two-Echelon Location Routing Problem}
\newacronym{2evrp}{2E-VRP}{Two-Echelon Vehicle Routing Problem}
\newacronym{2tcl}{2T-CL}{Two-Tier City Logistics}
\newacronym{2tssnd}{2T-SSND}{Two-Tier Scheduled Service Network Design}
\newacronym{sdp}{VRPSDPTW}{Vehicle Routing Problem with Simultaneous Delivery and Pickup and Time Windows}
\newacronym{lp}{LP}{Linear Programming}
\newacronym{bpc}{BPC}{Branch-Price-and-Cut}

%% file: CLIG_1-Intro.tex
% INTRODUCTION
The distribution of goods is a vital component of modern economies, supporting economic and social activities. However, urban freight transportation significantly contributes to negative externalities such as traffic congestion, noise, air pollution, and emissions, which adversely affect urban environments and public health \citep{savelsbergh201650th, demir2015selected}. This situation is expected to intensify as urbanization accelerates and the freight volume caused by growing e-commerce and local consumption increases \citep{EuropeanCommission2023}.

City logistics aims at designing and managing efficient freight transportation systems in urban areas so that the negative impacts are minimized, while supporting economic and social development and treating all stakeholders as part of an integrated system \citep{taniguchi2001CL, crainic2009models}. A well-established strategy to achieve this goal is the operation of consolidation-based two-tier transportation networks with heterogeneous and multi-modal vehicle fleets \citep{savelsbergh201650th, benmohamed2017InterconCL, boysenCLSurvey2020}. The two-tier structure is particularly beneficial if customer density is high and distances between depots and distribution areas are long \citep{fontaine2023smart}. In addition, it reduces transportation cost and easily integrates existing infrastructure, like public transport, as well as environmentally friendly transport modes, like cargo bikes and delivery on foot.

\subsection{Two-Tier City Logistics}

The concept of \textit{\gls{2tcl}} was first introduced by \citet{crainic2004advanced} as a location-allocation problem for satellite depots, laying the groundwork for research at all planning levels.
Since \gls{2tcl} covers a wide range of planning problems, the terminology in the literature is not consistent. While we adopt the notation of \citet{crainic2004advanced}, we give the corresponding alternatives wherever possible.
In \gls{2tcl}, a set of demands representing parcels with space- and time-related requirements must be shipped through a transportation network with two layers \citep{crainic2004advanced}. As shown in Figure \ref{fig:PhyNetwork}, \textit{\glspl{dc}}, also called city or urban distribution centers \citep{fontaine2021scheduled,dumez2023matheuristic} or external zones \citep{crainic2020planning}, are located easily accessible on the outskirts of the city. Inbound and outbound freight is handled and (de-)consolidated before or after its shipment in the inner city area, where customer locations are situated. \textit{Satellites} denote smaller depots or flexible locations that provide additional (de-)consolidation points and allow operators to switch to vehicle types compatible with the urban area. In the literature, satellites are also referred to as urban hubs \citep{benmohamed2017InterconCL}, intermediate depots \citep{zhao2018heterogeneous} or micro hubs \citep{boysen2023optimization}. Transportation on the \textit{first tier} occurs between the \glspl{dc} and the satellites, where large vehicles such as trucks or trams operate. The \textit{second tier} connects satellites and customer locations using smaller, environmentally friendly vehicles like cargo bikes.
\ifbool{TR}
{Since \gls{2tcl} covers a wide range of planning problems, the terminology in the literature is not consistent. While we adopt the notation of \citet{crainic2004advanced}, we give the corresponding alternatives wherever possible.}
{}
Tiers are also called levels or echelons and, referring to the mostly circular structure of urban areas, some studies use the terms outer and inner to describe the first and second transport layer, respectively, \citep[see, e.g.,][]{grangier2016adaptive, anderluh2021multi}. 

\myfigure[]
{Physical Network Representations of a \gls{2tcl} System}
{fig:PhyNetwork}
{tbp}
{
	\scalebox{1}
	{
		\input{physical_network}
	}
}

In many studies on \gls{2tcl} planning problems, the goal is the minimization of the total operating cost comprising, e.g., fixed cost for the operation of vehicles and variable cost for the transport of the demands \citep{crainic2009models, benmohamed2017InterconCL, zhao2018heterogeneous, dumez2023matheuristic}. These costs do not only cover monetary aspects but may also include environmental impacts and factors related to sustainability or quality of life. The consideration of negative externalities of shipments, e.g., expressed in the total waiting time of the vehicles or the total carbon emissions, receives increasing attention \citep[see][]{anderluh2021multi, lange2025TechRep2TCLSS}.
Depending on (i) the planning horizon, i.e., strategic, tactical or operational, (ii) the decision hierarchy, i.e., which decisions are taken before and after the current planning step, and (iii) the scope of control, i.e., planning an own fleet or integrating externally offered transportation services, a variety of planning problems arise in \gls{2tcl}.
%As a basis for this work, we would like to mention three frequently studied ones.
%Three of the most frequently studied ones are the following.
The following three basic problems are among the most frequently studied ones.
%In the \textit{\gls{2evrp}}, two independent fleets of vehicles operate on the two tiers and a route for each vehicle must be found so that all demands are satisfied. Thus, routing decisions are taken defining whether a certain vehicle travels between a certain pair of locations.
In the \textit{\gls{2evrp}}, two independent fleets of vehicles operate on the two tiers to satisfy a given set of demands. The task is to plan a route for each vehicle such that the overall cost is minimum.
The \textit{\gls{2elrp}} extends the \gls{2evrp} by additionally deciding on which \gls{dc} or satellite locations to open or close.
In contrast, \textit{\gls{2tssnd}} takes a set of transportation services with predefined routes and schedules as input for each tier.
The task is to select a set of services to operate and to determine a scheduled itinerary for each demand.
%These settings can also be applied to individual tiers and combined, and the planning problems can be extended by various characteristics, such as demand \glspl{tw}, different commodity flow directions, maximum numbers of vehicles or resource management, cost allocation, and inter-tier synchronization aspects \citep[see, e.g.,][]{grangier2016adaptive, zhao2018heterogeneous, fontaine2021scheduled, dumez2023matheuristic, guckel2024resourceCost, lange2025TechRep2TCLSS}.
These basic problem characteristics (routing vs.\ service network design) can also be applied to individual tiers and combined. Furthermore, the planning problems can be extended by various characteristics, such as demand \glspl{tw}, different commodity flow directions, maximum numbers of vehicles or resource management, cost allocation, and inter-tier synchronization aspects \citep[see, e.g.,][]{grangier2016adaptive, zhao2018heterogeneous, fontaine2021scheduled, dumez2023matheuristic, guckel2024resourceCost, lange2025TechRep2TCLSS}.

\subsection{Instance Generation in Transportation and Logistics}

Benchmark instances are supposed to facilitate scientific testing by meeting the requirements of diversity, realism, size, and extensibility \citep[see][]{vanhoucke2007nsplib}. The vehicle routing community features the longest tradition in using benchmark instances, where current studies are still based on the well-known instances published in \citet{solomon1987algorithms}. However, defining customer locations randomly on a grid and using Euclidean distances reflects real-world situations to a limited extend only. In addition, methodological and technical developments turned previously complex instances into easily solvable ones \citep[see][]{uchoa2017new} and new characteristics require the generation of particular parameters in each individual study. 
Other areas of operations research also tackled the need for benchmark libraries, see, e.g., 
the TSPLIB by \citet{reinelt1991TSPlib} for the traveling salesman problem, 
the PSPLIB by \citet{kolisch1997PSPLib} for project scheduling problems and the MIPLIB for mixed-integer programming (MIP) problems, which is regularly updated to meet the state-of-the-art \citep[see the latest version in][]{gleixner2021MIPLIB}.

Despite the extensive research work on \gls{2tcl} systems, the research community lacks standardized and real world-inspired benchmark data that is suitable for modern characteristics like multiple tiers and modes \citep[see][]{boysenCLSurvey2020}. Studies often rely on artificial data or specific cases, resulting in limited availability of data, applicability, consistency, and comparability of results. This complicates the evaluation of the performance of algorithms and research progress is slowed down by the repeated effort in data preparation and instance generation. 

\subsection{Contributions}

To overcome these issues, we contribute the following aspects:
\begin{itemize} 
	\item We summarize and classify necessary and optional characteristics of city logistics problems considered in the literature. The analysis highlights the complexity of the instance generation resulting from the variety of characteristics, their interdependencies, and feasibility issues. 
	\item We provide the CLLIB -- a first benchmark library for urban transportation systems with up to two tiers \citep{Lange2025CLLIB}. Based on real cities' topologies, it includes all necessary and many optional characteristics. Besides, an add-on library CLLIB-Analytics for statistical evaluations is set up.
	\item Since it is hardly possible to cover the variety of problem settings with a single library, we propose the general, well-documented, and publicly available instance generation framework {CLInstGen} \citep{Lange2025CLInstGen}. It is easily extendable and provides 15 control parameters to generate diverse and realistic instances.
	\item We demonstrate the applicability of the CLLIB to single- and two-tier settings through computational studies, showcase possibilities for adaptations, and provide insights on the effects of instance parameters on solvability by general MIP methods and solution structures.
\end{itemize}   

\subsection{Structure}

The remainder of the paper is organized as follows. Section \ref{sec:InstDes} summarizes relevant problem characteristics and their interdependencies. Section \ref{sec:InstGen} presents the instance generation framework, explaining requirements, the algorithmic scheme, the data used, and the structure of the CLLIB and its components. In Section \ref{sec:Study}, we apply the library to two different problems and analyze the effects of various instance characteristics. Section \ref{sec:Concl} concludes and outlines future research directions.

%% file: physical_network.tex
\tikzset{
	cust/.style={circle,fill=rptublgrey,text=white,font=\footnotesize,minimum size=14pt,inner sep=0pt},
	sat/.style={regular polygon,regular polygon sides=3,fill=rptured,text=white,font=\footnotesize,minimum size=24pt,inner sep=0pt},
	dc/.style={rectangle,fill=black,text=white,font=\footnotesize,minimum size=14pt,inner sep=0pt},
	tt1/.style={font=\footnotesize},
	tt2/.style={font=\footnotesize\bfseries},
	custG/.style={circle,fill=rptumblgrey,text=white,font=\footnotesize,minimum size=14pt,inner sep=0pt},
	satG/.style={regular polygon,regular polygon sides=3,fill=rptumblgrey,text=white,font=\footnotesize,minimum size=24pt,inner sep=0pt},
	dcG/.style={rectangle,fill=rptumblgrey,text=white,font=\footnotesize,minimum size=14pt,inner sep=0pt},
	tt1G/.style={text=rptumblgrey,font=\footnotesize},
	tt2G/.style={text=rptumblgrey,font=\footnotesize\bfseries},
}
\begin{tikzpicture}
	%inner city
	\node[color=gray] (in) at (1.5,0.75) {\parbox{30mm}{\centering\small inner city area}};
	%rural area
	\node[color=gray] (in) at (-5,0.75) {\parbox{14mm}{\centering\small rural area}};
	% legend
	\node at (7.85,-0.5) {\faIcon{warehouse}};
	\node[anchor=mid west] at (8.35,-0.5) {\small\gls{dc}};
	\node[rptured] at (7.85,-1.15) {\faIcon{exchange-alt}};
	\node[anchor=mid west] at (8.35,-1.15) {\small Satellite};
	\node[rptublgrey] at (7.85,-1.9) {\faIcon{home}/\faIcon{building}};
	\node[anchor=mid west] at (8.35,-1.9) {\parbox{15mm}{\small Customer\\[-1.2ex] Location}};
	%customers
	\node[rptublgrey] (c6) at (0,0.5) {\faIcon{home}};
	\node[rptublgrey] (c7) at (-3,-1.25) {\faIcon{home}};
	\node[rptublgrey] (c8) at (1,-1.5) {\faIcon{building}};
	\node[rptublgrey] (c9) at (2.5,0) {\faIcon{building}};
	\node[rptublgrey] (c10) at (4,-2.5) {\faIcon{home}};
	\node[rptublgrey] (c11) at (2,-3) {\faIcon{building}};
	\node[rptublgrey] (c12) at (-2,-1.75) {\faIcon{home}};
	\node[rptublgrey] (c13) at (4.5,-0.5) {\faIcon{home}};
	\node[rptublgrey] (c14) at (5.5,0) {\faIcon{home}};
	\node[rptublgrey] (c15) at (3.5,0.5) {\faIcon{building}};
%	\node[gray] at (-0.25,-1.5) {\faIcon{home}};
%	\node[gray] at (1.25,-0.75) {\large\faIcon{city}};
%	\node[gray] at (2.75,-1.25) {\large\faIcon{city}};
%	\node[gray] at (2.25,-1.5) {\faIcon{hospital}};
%	\node[gray] at (2.25,-3) {\large\faIcon{city}};
	%satellites
	\node[rptured] at (-3.5,0.25) {\faIcon{store-alt}};
	\node[rptured] (s1) at (-3,-0.1) {\faIcon{exchange-alt}};
	\node[rptured] (s2) at (0.5,-3) {\faIcon{exchange-alt}};
	\node[rptured] at (3.25,-0.55) {\faIcon{school}};
	\node[rptured] (s3) at (3,-1) {\faIcon{exchange-alt}};
	%distribution centers
	\node (d1) at (-5.5,-1) {\faIcon{warehouse}};
	\node (d2) at (6.5,-3) {\faIcon{warehouse}};
	% rural areas
	%left
	\node[gray] at (-5,-1.75) {\large\faIcon{industry}};
	\node[gray] at (-4,-1) {\large\faIcon{tree}};
	\node[gray] at (-3.75,-2.25) {\large\faIcon{tree}};
	\node[gray] at (-5.25,-0.25) {\large\faIcon{tree}};
	\node[gray] at (-5,0) {\large\faIcon{tree}};
	\node[gray] at (-4.5,-0.1) {\large\faIcon{tree}};
	% right
	\node[gray] at (6.25,-2.25) {\large\faIcon{industry}};
	\node[gray] at (5.85,-1) {\large\faIcon{tree}};
	\node[gray] at (6.15,-1.25) {\large\faIcon{tree}};
	%inner city edges (city freighters)
	\draw (s1) --node[tt1,above]{} (c6);
	\draw (s1) --node[tt1,above right]{} (c7);
	\draw (s1) --node[tt1,above right]{} (c12);
	\draw (s1) --node[tt1,above right]{} (c8);
	\draw (c6) --node[tt1,right]{} (c7);
	\draw (c6) to[out=-35,in=90]node[tt1,right]{} (c11);
	\draw (c6) --node[tt1,right]{} (c12);
	\draw (c7) --node[tt1,below]{} (c8);
	\draw (c7) --node[tt1,below]{} (c12);
	\draw (c6) --node[tt1,right]{} (c8);
	\draw (c6) --node[tt1,above]{} (c9);
	\draw (c8) --node[tt1,right]{} (c9);
	\draw (c8) --node[tt1,above]{} (c10);
	\draw (c8) --node[tt1,above]{} (c12);
	\draw (c8) --node[tt1,right]{} (c11);
	\draw (c9) to[out=-100,in=145]node[tt1,right]{} (c10);
	\draw (c9) --node[tt1,above]{} (s3);
	\draw (c9) --node[tt1,above]{} (c13);
	\draw (c9) --node[tt1,above]{} (c14);
	\draw (c9) --node[tt1,above]{} (c15);
	\draw (c10) --node[tt1,below]{} (s3);
	\draw (c10) --node[tt1,right]{} (c11);
	\draw (c10) --node[tt1,above]{} (c13);
	\draw (c10) --node[tt1,above]{} (c14);
	\draw (c10) --node[tt1,above]{} (c15);
	\draw (c8) --node[tt1,above left]{} (s2);
	\draw (c11) --node[tt1,above left]{} (s2);
	\draw (c12) --node[tt1,above left]{} (s2);
	\draw (c13) --node[tt1,above]{} (c14);
	\draw (c13) --node[tt1,above]{} (c15);
	\draw (c13) --node[tt1,above]{} (c11);
	\draw (c14) --node[tt1,above]{} (c15);
	\node at (3.35,-3) {\small\faIcon{plug} \faIcon{car-side}};
	\node at (1.5,-0.2) {\faIcon{biking}};
	\node at (-0.8,-1.15) {\faIcon{biking}};
	\node[rptured] at (0.5,-2.5) {\faIcon{gas-pump}};
	%external edges (urban vehicles)
	% road-based
	\draw[double, double distance=3pt] (d1) --node[tt1,left]{} (s1);
	\draw[dashed,thick] (d1) --node[tt1,left]{} (s1);
	\draw[double, double distance=3pt] (d1) --node[tt1,below]{} (s2);
	\draw[dashed,thick] (d1) --node[tt1,below]{} (s2);
	\draw[double, double distance=3pt] (s1) to node[tt1,below left]{} (s2);
	\draw[dashed,thick] (s1) to node[tt1,below left]{} (s2);
	\draw[double, double distance=3pt] (s2) to node[tt1,below right]{} (s3);
	\draw[dashed,thick] (s2) to node[tt1,below right]{} (s3);
	\node at (-2.75,-2.5) {\faIcon{truck}};
	% rail-based
	\draw[draw=black,double distance=3pt,line width=1pt] (d2) --node[tt1,right]{} (s3);
	\draw[draw=gray,densely dashed,line width=4pt] (d2) --node[tt1,right]{} (s3);
	\draw[draw=black,double distance=3pt,line width=1pt] (s1) to node[tt1,above]{} (s3);
	\draw[draw=gray,densely dashed,line width=4pt] (s1) to node[tt1,above]{} (s3);
	\node at (5.25,-1.75) {\faIcon{subway}};
	
\end{tikzpicture}

%% file: CLIG_3-InstDes.tex
% INSTANCE DESCRIPTION

Many problems in \gls{2tcl} share equivalent or related characteristics that are categorized and structured in the following sections. We characterize problems based on the facilities (Section \ref{Facilities}), customer demands (Section \ref{Customer_demands}), as well as first- and second-tier transportation (Section \ref{services}). We further derive the corresponding instance requirements implied by these problem settings, distinguishing between those that are essential, i.e., required in most problem settings, and those that are optional, i.e., considered only in some settings. Finally, we present an overview of problem characteristics used in the literature in Section \ref{Instances_literature}.

\subsection{Facilities}\label{Facilities}
Despite inconsistent terminology, most problems consider a basic distinction between first- and second-tier facilities. In \gls{2tcl}, the physical infrastructure comprises \glspl{dc} and satellite depots, each with specific roles and requirements.

\textbf{\glspl{dc}}: These centers are typically located on the city's outskirts to ensure easy access for freight transportation. Depending on the problem setting, these locations can either be predetermined (fixed) or represent potential sites where the decision to establish a \gls{dc} is part of the problem. In the latter case, additional data, such as opening cost, may be required in strategic problem settings. The essential requirement is the location of the \glspl{dc}, while capacity restrictions are optional and only considered in some problem settings.

\textit{necessary: location \quad optional: capacity, opening cost}

\textbf{Satellites}: Satellites are smaller depots that serve as connecting points between the first and second tier. These depots are distributed within the city to ensure proximity to customer zones. The essential requirement is the location of the satellites. In many problem settings, capacity constraints are also considered, for example limits on the number of operating vehicles \cite[e.g.,][]{crainic2016modeling, guckel2024resourceCost} or volume limitations for demand handling \cite[e.g.,][]{fontaine2021scheduled}. Some problem settings further distinguish between different satellite types, such as road-based and rail-based satellites, which differ in their infrastructure characteristics \citep{lange2025TechRep2TCLSS}.

\textit{necessary: location \quad optional: capacity, maximum waiting time, opening cost}

\subsection{Customer Locations and Demands}\label{Customer_demands}
Customers are located within the city and may be associated with attributes such as a maximum waiting time, i.e., the time vehicles are allowed to wait at the location. Customer locations are typically distributed across the city to reflect realistic patterns, e.g., according to population density.

Demands represent the delivery requirements that must be fulfilled. Each demand is linked to a customer location, which serves as its destination. A distinction can be made regarding the origin of demands: while some problem settings assume demands originate in external zones \citep[e.g.,][]{crainic2009models}, others assume that demands have an origin point and must first be transported to \glspl{dc}, thereby incurring additional cost \citep[][]{fontaine2021scheduled, guckel2024resourceCost}. 

Demands may also be subject to time-related constraints, such as a release date or \gls{tw} when they become available at a \gls{dc}, and a due date or \gls{tw} for arrival at the final destination. Furthermore, demands can be categorized as inbound or outbound, depending on whether the freight is entering or leaving the city \citep[e.g.,][]{fontaine2021scheduled}. Some problem settings also account for the allocation of demands to specific logistics service providers \citep[e.g.,][]{guckel2024resourceCost, crainic2020planning}. Additional cost components may arise, such as allocation cost for assigning demands to \glspl{dc} or satellite handling/usage cost \cite[][]{crainic2016modeling}.

\textit{necessary: location, volume \quad optional: time dimension, types (e.g., inbound/outbound), maximum waiting time}

\subsection{First and Second Tier Transportation}\label{services}
\textbf{First tier}: Depending on the problem setting, first-tier transportation is represented either as routing decisions based on distances or through service-based transportation using urban vehicle services. 

Services are pre-generated transport routes from \glspl{dc} to sequences of satellite depots, facilitating goods movement at the first tier of the \gls{2tcl} system. These services are particularly relevant in problem settings capturing scheduled transportation systems at the tactical level, whereas routing-based settings typically do not rely on such pre-generated services. One motivation for using services instead of routing decisions is that they better represent transportation modes such as trams, which operate according to predetermined timetables and visit a sequence of stops. 

Service characteristics include the sequence of satellite depots served, vehicle capacity, and associated fixed and variable costs, which may depend on the vehicle fleet. Time-related parameters include fixed departure and arrival times or time windows, as well as unloading times at each satellite. Some problem settings specify services as fixed routes, while others use service legs representing partial segments of a service. In the latter case, complete routes are constructed through the selection of service legs \citep[e.g.,][]{crainic2014service}. 

First-tier cost cover transportation from \glspl{dc} to satellites and are typically modeled either as distance-based routing cost or as service operating cost. In the latter case, each service has fixed and variable cost components, which may vary by vehicle type, while routing-based settings typically consider only distance-based cost.

\textit{necessary: distances and/or cost \quad optional: services, capacities, vehicle types}

\textbf{Second tier}: In most problem settings, second-tier transportation is represented either as routing decisions, approximations, or via pre-generated tours. Tours operate similarly to services but are associated with the second tier, carrying out final deliveries from satellite depots to customer locations. These tours are defined by the sequence of customer locations, fixed and variable costs, capacity constraints, and potential time-related requirements. 

System-wide constraints may limit the number of vehicles that can operate simultaneously across both tiers \cite[e.g.,][]{lange2025TechRep2TCLSS,perboli2011two}. 

Second-tier cost vary depending on the problem setting. In some approaches, they are approximated based on distances between satellites and customers or by incorporating satellite handling cost \citep[][]{fontaine2021scheduled}. In problem settings that explicitly model second-tier routing, cost are based on the distance traveled by vehicles \citep[][]{perboli2011two, guckel2024resourceCost}. Other approaches using pre-generated tours assign each tour a cost rate including fixed and variable components \citep{lange2025TechRep2TCLSS}.

\textit{necessary: distances and/or cost \quad optional: services, capacities, vehicle types}

\subsection{Problem Characteristics in the Literature}\label{Instances_literature}
In this section, we classify \gls{2tcl} publications based on key problem characteristics related to \glspl{dc}, satellites, demands, and first- and second-tier transportation. The overview covers a selection of works from service network design, vehicle routing, and location routing, and highlights which characteristics are considered in different problem settings. While not exhaustive, it provides a structured view of how these problems are defined in the literature. 

We also show that our instance generation framework can represent the requirements of these different problem settings. The classification is summarized in Table \ref{tab:rotated_table}.

\ifbool{TR}
{
	\afterpage{
		\begin{landscape}
			\begin{table}
				\caption{Classification of Literature}
				\label{tab:rotated_table}
				\scalebox{0.8}{
					\input{literature_table.tex}

				}
			\end{table}
		\end{landscape}
	}
}
{
	\mytable[] % <- no abbreviations here
	{Classification of Literature}
	{tab:rotated_table}
	{tbp}
	{
			\scalebox{0.69}{
				\input{literature_table.tex}
			}}
}

%% file: literature_table.tex
%\begin{tabular}{l|cc|cccc|ccccc|cccc|cccc}
\begin{tabular}{lccccccccccccccccccc}
	\toprule
	
	%& \multicolumn{2}{c|}{\textbf{DCs}} 
	%& \multicolumn{4}{c|}{\textbf{Satellites}} 
	%& \multicolumn{5}{c|}{\textbf{Customers / Demands}} 
	%& \multicolumn{4}{c|}{\textbf{First tier}} 
	%& \multicolumn{4}{c}{\textbf{Second tier}} \\
	& \multicolumn{2}{c}{\textbf{DCs}} 
	& \multicolumn{4}{c}{\textbf{Satellites}} 
	& \multicolumn{5}{c}{\textbf{Customers / Demands}} 
	& \multicolumn{4}{c}{\textbf{First Tier}} 
	& \multicolumn{4}{c}{\textbf{Second Tier}} \\
	\cmidrule(r){2-3} 
	\cmidrule(lr){4-7} 
	\cmidrule(lr){8-12}  
	\cmidrule(lr){13-16} 
	\cmidrule(lr){17-20}
	\textbf{Publication} 
	& Loc. & Cap. 
	& Loc. & Cap. & Wait. & Cost 
	& Loc. & Vol. & Time & Wait. & Type 
	& Dist. & Serv. & Cap. & Type 
	& Dist. & Serv. & Cap. & Type \\
	\midrule
	% ============================
	% Service Network Design
	% ============================
	\multicolumn{20}{l}{\textbf{Service Network Design}} \\
	\cite{crainic2009models} & \ding{51} & - & \ding{51} & \ding{51} & - & - & \ding{51} & \ding{51} & \ding{51} & - & - & - & \ding{51} & \ding{51} &  \ding{51}  & \ding{51} & - & - & -\\
	\cite{crainic2016modeling}  & \ding{51} & - & \ding{51} & \ding{51} & - & - & \ding{51} & \ding{51} & \ding{51} & - & - & - & \ding{51} & \ding{51} &  \ding{51} & \ding{51} & - & - & -\\
	\cite{fontaine2021scheduled} &  \ding{51} & \ding{51} & \ding{51} & \ding{51} & - & - & \ding{51} & \ding{51} & \ding{51} & - & \ding{51} &  - & \ding{51} & \ding{51} & \ding{51} & \ding{51} & - & - & -\\
	\cite{guckel2024MultiDay} &  \ding{51} & \ding{51} & \ding{51} & \ding{51} & - & - & \ding{51} & \ding{51} & \ding{51} & - & - &  - & \ding{51} & \ding{51} & \ding{51} & \ding{51} & - & - & -\\
	\cite{guckel2024resourceCost} &  \ding{51} & \ding{51} & \ding{51} & \ding{51} & - & - & \ding{51} & \ding{51} & \ding{51} & - & - &  - & \ding{51} & \ding{51} & \ding{51} & \ding{51} & - & \ding{51} & -\\
	\cite{lange2025TechRep2TCLSS} & \ding{51} & - & \ding{51} & \ding{51} & \ding{51} & - & \ding{51} & \ding{51} & \ding{51} & \ding{51} & \ding{51} & \ding{51} & \ding{51} & \ding{51} & - & \ding{51} & \ding{51} & \ding{51} & -\\
	\cite{masson2017optimization} & \ding{51} & -  & \ding{51} & \ding{51} & - & - & \ding{51} & \ding{51} & \ding{51} & - & - & - & \ding{51} & \ding{51} & - & \ding{51} & - & \ding{51} & -\\
	%\cmidrule(r){2-3} 
	%\cmidrule(lr){4-7} 
	%\cmidrule(lr){8-12}  
	%\cmidrule(lr){13-16} 
	%\cmidrule(lr){17-20}
	\cmidrule(){1-20}
	% ============================
	% Vehicle Routing
	% ============================
	\multicolumn{20}{l}{\textbf{Vehicle Routing}} \\
	\cite{anderluh2021multi} &  \ding{51} & - & \ding{51} & \ding{51} & \ding{51} & - & \ding{51} & \ding{51} & \ding{51} & - & - & \ding{51} & - & \ding{51} & - & \ding{51} & - & \ding{51} & -\\
	\cite{dellaert2019branch} &  \ding{51} & - &  \ding{51} & - &  \ding{51} & - &  \ding{51} &  \ding{51} &  \ding{51} & - & - &  \ding{51} & - &  \ding{51} & - &  \ding{51} & - &  \ding{51} & -\\
	\cite{dumez2023matheuristic} & \ding{51} & - & \ding{51} & \ding{51} & - & - & \ding{51} & \ding{51} & \ding{51} & - & \ding{51} & \ding{51} & - & \ding{51} & - & \ding{51} & - & \ding{51} & -\\
	\cite{grangier2016adaptive} &\ding{51}  & - & \ding{51} & \ding{51} & - & - & \ding{51} & \ding{51} & \ding{51} & - & - & \ding{51} & - & \ding{51} & - & \ding{51} & - & \ding{51} & -\\
	\cite{perboli2011two} & \ding{51} & \ding{51} & \ding{51} & \ding{51} & - & - & \ding{51} & \ding{51} & - & - & - & \ding{51} & - & \ding{51} &  - & \ding{51} & - & \ding{51} & -\\
	%\cmidrule(r){2-3} 
	%\cmidrule(lr){4-7} 
	%\cmidrule(lr){8-12}  
	%\cmidrule(lr){13-16} 
	%\cmidrule(lr){17-20}
	\cmidrule(){1-20}
	% ============================
	% Location Routing
	% ============================
	\multicolumn{20}{l}{\textbf{Location Routing}} \\
	\cite{mirhedayatian2021two} &   \ding{51} & - &  \ding{51}  & - & - &  \ding{51}  &  \ding{51}  &  \ding{51}  &  \ding{51}  & - &  \ding{51}  & - & - & - &  \ding{51}  & - &  \ding{51}  &  \ding{51} \\
	\cite{yildiz2023exact} & \ding{51} & - & \ding{51} & \ding{51} & - & \ding{51} & \ding{51} & \ding{51} & - & - & \ding{51} & - & \ding{51} & - & \ding{51} & - & \ding{51} & -\\
	\cite{zhao2018heterogeneous} & \ding{51} & - & \ding{51} & \ding{51} & - & \ding{51} & \ding{51} & \ding{51} & - & - & - & \ding{51} & - & - & \ding{51} & - & \ding{51} & \ding{51}\\
	%\cmidrule(r){2-3} 
	%\cmidrule(lr){4-7} 
	%\cmidrule(lr){8-12}  
	%\cmidrule(lr){13-16} 
	%\cmidrule(lr){17-20}
	\cmidrule(){1-20}
	% ============================
	% Instance Generation Framework
	% ============================
	\textbf{Instance Generation Framework} & \ding{51} & - & \ding{51} & \ding{51} & \ding{51} & - & \ding{51} & \ding{51} & \ding{51} & \ding{51} & \ding{51} & \ding{51} & \ding{51} & \ding{51} &  \ding{51} & \ding{51} & \ding{51} & \ding{51} &  \ding{51}\\
	\bottomrule
\end{tabular}

%% file: CLIG_4-InstGen.tex
% INSTANCE GENERATION PROCEDURE

This section provides detailed explanations of version v1 of the instance generation framework CLInstGen \citep{Lange2025CLInstGen}, starting with the instance characteristics covered (Section \ref{subsec:instGenCharac}) and the city data as well as real-world system parameters of vehicles and demands used (Sections \ref{subsec:cityData} and \ref{subsec:commonData}). The control parameters and the workflow of the instance generator are described in Sections \ref{subsec:instGenParam} and \ref{subsec:generator}, respectively. Due to a lack of empirical data for time-related demand characteristics and service properties, they are generated rule-based with reasonable interrelations as proposed in Section \ref{subsec:instGenRules}. The CLLIB is presented in Section \ref{subsec:CLLIB} and Section \ref{subsec:TightVehCap} provides an instance adaption scheme for research applications.

\subsection{Instance Characteristics}
\label{subsec:instGenCharac}

The data describing each benchmark instance is selected based on the necessary and optional characteristics of the \gls{2tcl} system identified in the previous section. As summarized in Table \ref{tab:rotated_table}, all necessary location characteristics are included (for \glspl{dc}, satellites, and customers). In line with the literature, detailed properties of \glspl{dc} are neglected, while characteristics of satellites are broadly represented. In particular, we provide satellite capacity restrictions measured by numbers of vehicles, size of the parking area, as well as maximum demand volume per time period and corresponding conversion rules. Maximum waiting times for customer locations are included to support the applicability with city logistics goals. Each demand is described by an origin and a destination determining its type, a volume, a pickup and/or delivery \gls{tw}, and a handover as well as availability time, if applicable.
The two-tiered network is outlined by homogeneous fleets of capacitated vehicles, related distances, and by transportation services with routes, fixed cost, and start \glspl{tw}.

\subsection{City Data}
\label{subsec:cityData}

The instance generation framework is based on real urban topologies. Every city's data includes the locations of \glspl{dc}, satellites, and customers. As described above, we assume that no direct transports exist between \glspl{dc} and customers. Satellite characteristics must be provided based on real locations like supermarkets, gas stations, tram stations, or other dedicated sites. Satellite capacity can be measured in terms of number of vehicles, size of the parking area, or demand volume per time period. We assume a proportionate behavior of these measures,
\ifbool{TR}
{where an urban vehicle requires $12 ~m^2$ and a city freighter $5 ~m^2$ for parking. To relate the maximum demand volume per time period with the maximum number of vehicles, the following rule is applied:
	$$\text{max. demand volume per time period} = 600,000 ~cm^3 \cdot \log_2(\text{total max. number of vehicles}),$$ 
	where $600,000 ~cm^3$ refers to two demands of the largest existing size. }
{see \cite{Lange2025CLInstGen} for details. }

Road- and rail-based satellites differ in a few characteristics. For road-based satellites, a region must be specified, which describes its geographic position in the urban area. When instances with a reduced number of satellites are generated, the generation scheme assures that they are uniformly distributed among regions. Rail-based satellites feature a maximum frequency of visits per hour, which is considered in the generation of the first-tier services. When a public transport system is used for instance generation, a set of lines must additionally be provided.

With the framework, data of three urban regions in Germany is provided \citep{Lange2025CLInstGen}. The numbers of \glspl{dc} and satellites are selected in line with the literature \citep[see][]{grangier2016adaptive, benmohamed2017InterconCL, fontaine2021scheduled, li2023sharedTransport}. The district Oststadt of Karlsruhe features about 21,000 inhabitants, an area of 5 km$^2$, 2 \glspl{dc}, 10 road-based and 4 rail-based satellites, as well as public transport by tram lines. Ingolstadt, a city with about 141,000 inhabitants and 133 km$^2$, has 3 \glspl{dc} and 11 road-based satellites. The city of Munich features about 1.5 million inhabitants in an area of 311 km$^2$ and comes with 4 \glspl{dc}, 20 road-based and 10 rail-based satellites. Maps of the urban areas showing the locations of \glspl{dc} and satellites can be found in \ifbool{TR}{\ref{app:Maps}}{Appendix \ref{app:Maps}}.

\subsection{Common Real-World System Parameters} 
\label{subsec:commonData}

Common parameters used for demand sizes, handover times, and vehicles capacities are based on package and vehicle information available for \glspl{lsp} operating in Europe, like DHL, Hermes and UPS, and reports on pilot projects on the integration of public transport into distribution systems. Demands feature five size categories with handover times and volumes as given in Table \ref{tab:demandSize}. If no empirical demand size data is provided with the customer locations, the given default size frequencies are used. Four vehicle categories are identified with average capacities as follows: cargo bike ($2~m^3$), electric car ($4~m^3$), small truck ($11~m^3$), and tram ($21~m^3$), where the first two are interpreted as city freighters and the last two act as urban vehicles. 

Note that this relation of demand sizes and vehicle capacities results in large amounts of packages that can be transported by a single vehicle. Based on the given demand size distribution, a cargo bike fits about 70 and a truck may load more than 200 parcels. While this mirrors real-world scenarios \citep[see, e.g.,][]{DHL2019FactsScooter}, instances with 100 demands and more are typically challenging in many problem settings in the city logistics literature. Thus, in instances of tractable size, the vehicle capacities are expected to never be a limiting factor. We address this issue by proposing an appropriate scaling rule in Section \ref{subsec:TightVehCap}.

\subsection{Instance Generation Parameters}
\label{subsec:instGenParam}

The instance generation is controlled by the following 15 parameters, with several significant dependencies existing between them.

\begin{enumerate}
	\item \textit{City:} defines the urban region considered
	\item \textit{Number of customer locations:} defines the number of customer locations used as origins or destinations of demands
	\item \textit{Number of demands:} defines the number of demands; can be greater than or equal to the number of customers
	\item \textit{Inner-city access \gls{tw}:} bounds the allocation of demand \glspl{tw}, effects planning time horizon, demand availability times and starting \glspl{tw} of services
	\item \textit{Instance density:} defines the density of the demand \gls{tw} allocation over time and effects the planning time horizon
	\item \textit{Satellite type:} defines whether road- or rail-based satellite locations are used; defines the type of urban vehicles; effects the first-tier service generation scheme
	\item \textit{Amount of satellites:} defines the amount of satellites: as given, halved, or minimum; can only be halved or minimum for road-based satellites; effects the number of first-tier services
	\item \textit{Second-tier vehicle type:} defines the type of city freighters used; effects the second-tier distances given in travel time 
	\item \textit{Demand type distribution:} defines the portions of e2c, c2e, and c2c demands 
	\item \textit{Demand size data source:} defines whether demand sizes are taken as input with the customer locations or if the default distribution of demand sizes is used
	\item \textit{Demand \gls{tw} widths distribution:} defines the distribution of short, medium, and long \glspl{tw} 
	\item \textit{Customer maximum waiting time distribution:} defines the distribution of no, short, and long waiting time allowed at customer locations
	\item \textit{Indicator for service generation:} defines whether transportation services are generated
	\item \textit{Indicator for generation of feasibility-assuring services:} defines whether feasibility-assuring services are generated
	\item \textit{Number of instance replicates:} defines the number of samples drawn whenever a distribution is used to determine parameter values
\end{enumerate}

\textit{Cities} of different size and shape can be used, which particularly effects the relation of average travel times between locations and the inner-city access time window. The \textit{number of customer locations} and the \textit{number of demands} can be controlled separately so that scenarios with single-family homes and apartment buildings are covered. The \textit{inner-city access \gls{tw}} can be varied freely so that instances of different demand density can be created. An appropriate combination of \textit{number of demands} and \textit{inner-city access \gls{tw}} must be chosen to avoid the generation of infeasible instances, where the distances implied by the urban region are also relevant. The demand \glspl{tw} are allocated randomly over the \textit{inner-city access \gls{tw}}, while demand availability times and travel times are considered if necessary. A planning time horizon, in which vehicles are allowed to travel to support instance feasibility, is defined by adding identical lead and follow-up times to the \textit{inner-city access \gls{tw}}. Nonetheless, capacity restrictions for vehicles as well as satellites and the number of services may still cause infeasibility. Alternative to an explicit \textit{inner-city access \gls{tw}}, the \textit{instance density} can be defined as low, medium, or high. The choice of the \textit{satellite type} and first-tier services can significantly effect the structure of the instance and the generation scheme. For road-based satellites, the \textit{amount of satellites} can be reduced to a halved or minimum amount starting from all provided (potential) satellites as the maximum set. Since the resulting number of satellites restricts the set of constructible first-tier service routes, the number of services decreases with the amount of satellites. The \textit{second-tier vehicle type} defines whether cargo bikes or cars are used in the dense urban area with corresponding travel times between satellites and customer locations.

The \textit{demand type distribution} controls whether and how many e2c, c2e, and c2c demands are included.
Since every portion can be set to zero, extreme demand type scenarios can be created just as average and real-world cases.
Note that in choosing e2c demands only, the standard case of planning the commodity flow from external zones into the city is implemented. There are two alternatives regarding the \textit{demand size data source}.
Either a demand size is given with every customer location from empirical data or the size of every demand is randomly generated according a distribution (see Section \ref{subsec:commonData}).
The \textit{demand \gls{tw} widths distribution} determines how many \glspl{tw} with short, medium, and long width exist.
Note that the long \gls{tw} option implies no restriction so that customers without \glspl{tw} can be integrated.
Similarly, the \textit{customer maximum waiting time distribution} defines how many customer locations allow no waiting, a short waiting time, or infinite waiting. 

The \textit{indicator for service generation} and the \textit{indicator for generation of feasibility-assuring services} define, respectively, whether transportation services are generated and whether the special subset of feasibility-assuring services (see Section \ref{subsec:instGenRules}) is included. This provides options to reduce the computational effort of the instance generation, when services are not required for routing-based problems or time-related infeasibility of instances cannot occur.

The \textit{number of instance replicates} defines the number of different instances generated with the same parameter setting.
This is particularly useful, when statistical reliability is required for analyses of effects of changes in single instance parameters.
\textit{Number of instance replicates} many samples are drawn for the selection of customer locations, reduced amounts of satellites, demand sizes taken from the default distribution, demand \gls{tw} widths, and maximum waiting times at customer locations.
Every replicate is combined with every other replicate. Thus, the number of instances increases enormously with the \textit{number of instance replicates}. 

Setting these parameters, in particular, considering different amounts of satellites, generating predefined services, and varying demand type as well as demand \gls{tw} widths distributions, gives rise to applications in the areas of vehicle routing, location routing, and service network design. Furthermore, the use of proper sampling and instance replicates provides the opportunity of generating sets of benchmark instance with a high reliability in statistical analyses of results.

\subsection{Instance Generator}
\label{subsec:generator}

The dependencies between parameters reported in Section \ref{subsec:instGenParam} imply precedence relations between instance generation steps. Together with the data describing the urban region, this is the main driver of complexity in the instance generation. Figure \ref{fig:InstGenFramew} shows the workflow of the instance generator, starting from data and parameter input depicted at the top, indicating the main generation steps and their dependencies in the middle, and ending with instance data output at the bottom. The first three steps use city data and instance parameters to read and select \gls{dc}, satellite, and customer locations. The region characteristic of the satellites is considered if a reduced amount of satellites is required. For a minimum amount of satellites, one satellite per existing region is randomly chosen. If all satellites feature the same region index, instances with one satellite are created. To implement half the amount of satellites, for each region, the number of satellites is divided by two and rounded down, and the resulting amount is randomly selected. The satellite locations are used to determine the distance matrix and act as input for the demand and service generation.

The framework uses \texttt{openrouteservice by HeiGIT} to obtain distance matrices, while distances are measured as travel times in minutes. The distances are used to generate demands as well as services and the inner-city access \gls{tw}, if an instance density parameter is given. To obtain an inner-city access \gls{tw} from an instance density parameter, an amount of travel hours is determined by summing up the longest direct travel distances from each location to any other location and rounding up to an integer amount. This amount of travel hours equals the inner-city access \gls{tw} for medium instance density and it is halved and doubled for high and low density instances, respectively. To realize transport connections beyond the urban region and support instance feasibility, vehicles are allowed to travel before and after the inner-city access \gls{tw}.
\ifbool{TR}
{Therefore, a planning time horizon is defined by adding identical lead and follow-up times. The length of the lead and follow-up times is determined based on the average travel times between all \glspl{dc} and satellites and between all satellites and customer locations. These average travel times are summed and rounded to the next greater multiple of 30 minutes.}
{Therefore, a planning time horizon is defined by adding identical lead and follow-up times, see \cite{Lange2025CLInstGen} for details. }

Since the inner-city access \gls{tw} defines the maximum demand \gls{tw} widths and the planning time horizon is used to set the maximum waiting time at customer locations, the maximum service start \gls{tw} and the demand availability times, they are required in three following steps. Maximum waiting times at customer locations are generated according to the given distribution. The short maximum waiting time is determined as $10/3$ times the maximum existing handover time, equaling 10 minutes with default values. The generation of the demands requires the most input parameters and data as it defines various location- and time-related characteristics. In the generation of services, the distances are used together with the planning time horizon and the locations to avoid unreasonably long service routes. Details on the rules are given in the following section.

\myfigure[]
{Workflow of the Instance Generator}
{fig:InstGenFramew}
{tbp}
{
\begin{tikzpicture}
		% INPUT
		\node[rectangle, fill=white!93!black] (citydata) at (0,-0.6) {\footnotesize City Data};
		\node[rectangle, fill=white!80!black] (satdata) at (1,-2.1) {\parbox{15mm}{\centering\footnotesize Sat. Type and Amount}};
		\node[rectangle, fill=white!80!black] (numcust) at (3.2,-2.1) {\parbox{17mm}{\centering\footnotesize No. of Customers}};
		\node[rectangle, fill=white!80!black] (stvehicle) at (5.4,-2.1) {\parbox{16mm}{\centering\footnotesize Sec.-Tier Vehicle Type}};
		\node[rectangle, fill=white!80!black] (th) at (7.4,-1.85) {\parbox{15mm}{\centering\footnotesize I.-C. Access \gls{tw}/ Inst. Dens.}};
		\node[rectangle, fill=white!80!black] (custwait) at (9.5,-2.1) {\parbox{18mm}{\centering\footnotesize Cust. Max. Wait. Distr.}};
		\node[rectangle, fill=white!80!black] (deminput) at (12.4,-2.1) {\parbox{30mm}{\centering\footnotesize Demand Amount, Type and TW Widths Distr.}};
		\node[rectangle, fill=white!93!black] (demsize) at (13,-0.6) {\footnotesize Demand Size (Distr.)};
		% FUNCTIONS
		\node[rectangle, fill=white!80!green] (readCDC) at (-0.2,-4.3) {\parbox{10mm}{\centering\footnotesize read CDC loc.}};
		\node[rectangle, fill=white!80!green] (genSat) at (1,-5.8) {\parbox{15mm}{\centering\footnotesize generate sat. loc.}};
		\node[rectangle, fill=white!80!green] (genCust) at (3,-5.8) {\parbox{15mm}{\centering\footnotesize generate cust. loc.}};
		\node[rectangle, fill=white!80!green] (genDist) at (5.4,-6.5) {\parbox{15mm}{\centering\footnotesize generate distances}};
		\node[rectangle, fill=white!80!green] (genTH) at (7.8,-6.25) {\parbox{22mm}{\centering\footnotesize read/generate access/horizon parameters}};
		\node[rectangle, fill=white!80!green] (genWait) at (10.6,-5.8) {\parbox{22mm}{\centering\footnotesize generate cust. max. wait.}};
		\node[rectangle, fill=white!80!green] (genDem) at (12.4,-8.6) {\parbox{15mm}{\centering\footnotesize generate demands}};
		\node[rectangle, fill=white!80!green] (genServ) at (13.8,-6.85) {\parbox{13mm}{\centering\footnotesize \textcolor{white!80!green}{text} generate services \textcolor{white!80!green}{text}}};		
		% OUTPUT
		\node[rectangle, fill=white!80!blue] (locout) at (1.5,-9.8) {\parbox{40mm}{\centering\footnotesize Locations and Properties}};
		\node[rectangle, fill=white!80!blue] (distout) at (5.4,-9.6) {\parbox{15mm}{\centering\footnotesize Distance Matrix}};
		\node[rectangle, fill=white!80!blue] (thout) at (8,-9.6) {\parbox{27mm}{\centering\footnotesize I.-C. Access \gls{tw}, Plan. Time Horizon}};
		\node[rectangle, fill=white!80!blue] (demout) at (11.8,-9.8) {\parbox{15mm}{\centering\footnotesize Demands}};
		\node[rectangle, fill=white!80!blue] (servout) at (13.8,-9.8) {\parbox{15mm}{\centering\footnotesize Services}};		
		% INPUT ARCS
		\draw[-latex,white!70!black] (-0.6,-0.9) -- (-0.6,-3.5);
		\draw[-latex,white!70!black] (-0.4,-0.9) -- (-0.4,-3.3) -- (13.8,-3.3) -- (genServ);
		\draw[white,fill=white] (0.95,-3.25) rectangle (1.05,-3.35);
		\draw[white,fill=white] (2.93,-3.25) rectangle (3.03,-3.35);
		\draw[white,fill=white] (5.35,-3.25) rectangle (5.45,-3.35);
		\draw[white,fill=white] (7.35,-3.25) rectangle (7.45,-3.35);
		\draw[white,fill=white] (10.55,-3.25) rectangle (10.65,-3.35);
		\draw[white,fill=white] (12.35,-3.25) rectangle (12.45,-3.35);
		\draw[-latex,white!70!black] (-0.2,-0.9) -- (-0.2,-3.1) -- (1,-3.1) -- (genSat);
		\draw[-latex,white!70!black] (satdata) -- (genSat);
		\draw[-latex,white!70!black] (0,-0.9) -- (0,-1.1) -- (2,-1.1) -- (2,-3.1) -- (3,-3.1) -- (genCust);
		\draw[-latex,white!70!black] (3,-2.6) -- (genCust);
		\draw[-latex,white!70!black] (9.5,-2.6) -- (9.5,-3.1) -- (10.6,-3.1) -- (genWait);
		\draw[-latex,white!70!black] (3.2,-2.6) -- (3.2,-3.1) -- (10.6,-3.1) -- (genWait);
		\draw[white,fill=white] (5.35,-3.05) rectangle (5.45,-3.15);
		\draw[white,fill=white] (7.35,-3.05) rectangle (7.45,-3.15);
		\draw[-latex,white!70!black] (stvehicle) -- (genDist);
		\draw[-latex,white!70!black] (th) -- (7.4,-5.45);
		\draw[-latex,white!70!black] (deminput) -- (genDem);
		\draw[-latex,white!70!black] (14.3,-0.9) -- (14.3,-3.1) -- (12.4,-3.1) -- (genDem);		
		% OUTPUT ARCS
		\draw[-latex, dashed,white!70!black] (readCDC) -- (-0.2,-9.5);
		\draw[-latex, dashed,white!70!black] (genSat) -- (1,-9.5);
		\draw[-latex, dashed,white!70!black] (genCust) -- (3,-9.5);
		\draw[-latex, dashed,white!70!black] (genDist) -- (distout);
		\draw[-latex, dashed,white!70!black] (7.4,-7.1) -- (7.4,-9); % genTH -- th
		\draw[-latex, dashed,white!70!black] (genServ) -- (servout);
		\draw[white,fill=white] (5.3,-7.7) rectangle (5.5,-7.9);
		\draw[white,fill=white] (7.3,-7.7) rectangle (7.5,-7.9);
		\draw[-latex, dashed,white!70!black] (genWait) -- (10.6,-7.8) -- (3.4,-7.8) -- (3.4,-9.5);
		\draw[-latex, dashed,white!70!black] (genDem) -- (12.4,-9.55);	
		% INNER ARCS
		\draw[-latex] (readCDC) -- (4.25,-4.3) -- (4.25,-6.3) -- (4.5,-6.3);
		\draw[-latex] (readCDC) -- (12.7,-4.3) -- (12.7,-6.3) -- (13,-6.3);
		\draw[-latex] (genSat) -- (2,-5.8) -- (2,-6.7) -- (4.5,-6.7);
		\draw[-latex] (genSat) -- (2,-5.8) -- (2,-7.6) -- (13,-7.6); % genServ
		\draw[white,fill=white] (4,-6.65) rectangle (4.1,-6.75);
		\draw[white,fill=white] (4,-7.55) rectangle (4.1,-7.65);
		\draw[white,fill=white] (6.4,-7.55) rectangle (6.5,-7.65);
		\draw[-latex] (genCust) -- (4.05,-5.8) -- (4.05,-8.8) -- (11.5,-8.8); % genDem
		\draw[-latex] (genCust) -- (4.05,-5.8) -- (4.05,-7.4) -- (13,-7.4); % genServ
		\draw[white,fill=white] (6.4,-7.35) rectangle (6.5,-7.45);
		\draw[white,fill=white] (9.2,-7.35) rectangle (9.3,-7.45);
		\draw[white,fill=white] (9.2,-7.55) rectangle (9.3,-7.65);
		\draw[-latex] (genCust) -- (4.05,-5.8) -- (4.05,-6.5) -- (genDist);
		\draw[-latex] (genDist) -- (6.6,-6.5);
		\draw[-latex] (6.3,-6.7) -- (6.45,-6.7) -- (6.45,-7.2) -- (13,-7.2); % genDist -- genServ
		\draw[-latex] (6.3,-6.7) -- (6.45,-6.7) -- (6.45,-8.6) -- (11.5,-8.6); % genDist -- genDem
		\draw[white,fill=white] (9.2,-7.15) rectangle (9.3,-7.25);
		\draw[-latex] (9.05,-6.5) -- (13,-6.5); % genTH -- genServ
		\draw[-latex] (9.05,-6) -- (9.35,-6); % genTH -- genWait
		\draw[-latex] (9.05,-6.7) -- (9.25,-6.7) -- (9.25,-8.4) -- (11.5,-8.4); % genTH -- genDem		
	\end{tikzpicture}
}

\subsection{Rule-Based Generation of Demand and Service Characteristics}
\label{subsec:instGenRules}

\subsubsection{Demand Characteristics}

The framework aims at allocating demands uniformly over customer locations. Every customer location acts as a demand origin or destination at least once. Regarding the demand sizes, the common data described in Section \ref{subsec:commonData} is applied.

\mytable[]
{Demand Sizes, Handover Times and Relative Frequencies}
{tab:demandSize}
{tbp}
{
	\begin{tabular}{lrrrrr}
		\toprule
		\textbf{Demand Size} & \textbf{XS} & \textbf{S} & \textbf{M} & \textbf{L} & \textbf{XL} \\
		\midrule
		Volume (cm$^3$) & 4,800 & 10,000 & 27,000 & 100,000 & 300,000 \\
		\midrule
		Relative Frequency & 0.1 & 0.2 & 0.4 & 0.2 & 0.1 \\
		\midrule 
		Handover Time (min) & 1 & 1 & 2 & 2 & 3 \\
		\bottomrule
	\end{tabular}
}

Depending on the demand type, one or two \glspl{tw} are set for each demand. The widths are generated independently according to the \gls{tw} widths distribution as follows:
\begin{itemize}
	\item \textit{short} \glspl{tw} with a width of five times the maximum existing handover time (15 minutes with default values),
	\item \textit{medium} \glspl{tw} with a width of ten times the maximum existing handover time (30 minutes with default values), and
	\item \textit{long} \glspl{tw} with a width equal to the inner-city access \gls{tw}.
\end{itemize}
The demand \glspl{tw} are allocated randomly over the inner-city access \gls{tw}, while feasibility is assured by considering the availability times for e2c demands and the travel time for c2c demands. 

The \textit{availability time} of an e2c demand is randomly chosen according to the following rule: With a probability of 0.5, the availability time is zero. Otherwise, the availability time is taken from a triangular distribution with lower bound and mode equal to zero and upper bound determined as the minimum of (i) the middle of the planning time horizon and (ii) planning time horizon minus medium demand \gls{tw} width, average travel time between any \gls{dc} and satellite, and average travel time between any customer location and satellite.

To facilitate adaptability and extensibility of the instances, artificial inter-regional travel times between an e2c demand origin or a c2e demand destination and all \glspl{dc} are determined. The origin and destination locations are randomly generated within a 200 km radius around the city center coordinates of the given urban region.

\subsubsection{First-Tier Service Generation}

First-tier services can be road-based or rail-based.
The existence of tracks in the urban region is not necessary. Rail-based services are supposed to rely on a given network of public transport lines. The general goal of the generation rules is to assure sufficiently many inbound and outbound transportation options for every satellite from a time, capacity, and connectivity perspective. This is pursued by generating two subsets of services, namely, \textit{feasibility-assuring} and \textit{random} services. While feasibility-assuring services connect each satellite inbound and outbound to the nearest \gls{dc} with maximum start time flexibility, random services visit sequences of satellites and feature a tightly bounded start \gls{tw}. The service fixed cost are determined as 
$$\text{fixed cost } = \text{ total travel time } \cdot \text{ start time flexibility multiplier.}$$ 
As default values, we apply multipliers of one and two for start \gls{tw} widths less than and equal to half of the planning time horizon, respectively, and a multiplier of three for maximum flexibility with start \gls{tw} width equal to the planning time horizon.

\textit{Road-based services} can be cyclic or acyclic, but visit at least one \gls{dc}. The subset of feasibility-assuring services consists of one trivial inbound and outbound service for each satellite connecting it to the nearest \gls{dc} and featuring a start \gls{tw} width equal to the planning time horizon. The randomly generated services visit between one and four satellites, since longer routes are not considered reasonable from a practical perspective.
\ifbool{TR}
{For each number of visited satellites, different sequences are randomly generated, while sequences with particularly long total travel time are neglected.}
{} 
For each satellite, at least one sequence including it is chosen. Every chosen sequence is completed at the start and the end with the nearest \gls{dc}. The start \gls{tw} width multiplier and the number of start options are randomly determined according to the categories shown in Table \ref{tab:ServiceCategories}, where all categories are equally likely to appear. The width of the start \gls{tw} is determined as the product of the multiplier and the length of the planning time horizon. The start \glspl{tw} of services featuring the same sequence of satellites are generated non-overlapping and distributed over the planning time horizon. Note that for each start option, an individual service is generated so that the total number of random services may vary for each instance.

\mytable[]
{Overview of Categories of Road-Based Services and their Generation Properties}
{tab:ServiceCategories}
{tbp}
{
	\begin{tabular}{lccccc}
		\toprule
		\textbf{Category} & \textbf{XS} & \textbf{S} & \textbf{M} & \textbf{L} & \textbf{XL}   \\
		\midrule
		Start \gls{tw} width multiplier & $1/15$ & $1/10$ & $1/5$ & $1/2$ & $1$ \\
		Number of starts & $[1, 7]$ & $[1, 5]$ & 2 & 1 & 1 \\
		Flexibility multiplier & 1 & 1 & 1 & 2 & 3 \\
		\bottomrule
	\end{tabular}
}

\textit{Rail-based services} are generated based on existing public transport lines and properties of the tram or bus stations as satellites.
Note that open-source real-world data on public transport systems exists for more and more cities \citep[see, e.g.,][]{kujala2018PTdata}.
In the city data provided with the instance generation framework, rail-based service routes can be cyclic or acyclic, but all lines are connected to at least one \gls{dc}. This enables the assumption of forbidding demand handovers between first-tier services and represents a basic requirement for the service generation to work. Similar to the generation rules for road-based services, feasibility-assuring and random rail-based services are created. For every satellite, two trivial services connecting it to the nearest \gls{dc} in terms of line routes are set up. To generate the random services, one visiting line is randomly chosen for each satellite. For each chosen line, the maximum frequency per hour is determined based on the allowed maximum frequencies of the satellites involved. The maximum number of services following the given line is created from the maximum frequency per hour and the length of the planning time horizon. The default start \gls{tw} width for all random rail-based services is 5 minutes, which supposedly allows them to be fitted into an existing schedule. The final number of random services is instance-individual due to the random line selections and the maximum frequencies of the satellites.

\subsubsection{Second-Tier Service Generation}

We assume that second-tier services are operated by an own fleet of city freighters or highly flexible providers like bike couriers with a start \gls{tw} width equal to the planning time horizon. Thus, time-related infeasibility may mainly be caused by tight demand \glspl{tw} of commodities shipped by the same service. Therefore, the service generation rule must provide sufficiently many options to reach a customer location and avoid unrealistically long service routes. Second-tier service routes can start and end at the same or at different satellites. Their fixed cost are always equal to the total travel time of the service. 

Similar to the first-tier service generation scheme, \textit{feasibility-assuring} and \textit{random} second-tier services are created. The generation rule for feasibility-assuring services depends on the demand type. For every e2c and c2e demand, one trivial service is generated, which connects the corresponding customer location with the nearest satellite on a cyclic route. For every c2c demand, a trivial service directly connecting the pickup and delivery customer locations is generated, while it starts and ends at the nearest satellite. 
\ifbool{TR}
{
	The generation of random second-tier services starts with determining a \textit{reasonable service length} in terms of number of customer locations visited as
	$$\left\lceil\min\left\{2 \cdot \frac{\text{number of customers}}{\text{number of satellites}}~;~ \text{number of customers}\right\}\right\rceil.$$
	The generation scheme constructs as many random second-tier services as there exist customer locations. For each service, a length is drawn from a triangular distribution with lower bound equal to two, upper bound equal to the number of customers, and mode equal to the reasonable service length. A random sequence of customer locations of the given length is determined and the nearest satellites of the first and last customer location are added, respectively. Finally, the set of second-tier services is of cardinality number of demands plus number of customer locations.
}
{
	The generation of random second-tier services is based on a \textit{reasonable service length} that is determined in terms of number of customer locations visited by
	$$\left\lceil\min\left\{2 \cdot \frac{\text{number of customers}}{\text{number of satellites}}~;~ \text{number of customers}\right\}\right\rceil.$$
	The generation scheme constructs as many random second-tier services as there exist customer locations. A random sequence of customer locations of a selected length is determined and the nearest satellites of the first and last customer location are added, respectively. Finally, the set of second-tier services is of cardinality number of demands plus number of customer locations.
}

\subsection{CLLIB Benchmark Instances}
\label{subsec:CLLIB}

\subsubsection{Instance Parameter Values}

\mytable[]
{Summary of Instance Parameter Values applied in CLLIB}
{tab:InstParamValues}
{tbp}
{
	\ifbool{TR}
	{
		\scalebox{0.95}{\input{CLLIB_instance_parameter_table.tex}}
	}
	{
		\input{CLLIB_instance_parameter_table.tex}
	}
}

Table \ref{tab:InstParamValues} summarizes the instance parameter values used in the generation of version v1 of the CLLIB. The table differentiates the three urban areas, since the available data effects some parameters' values. We define the number of customer locations equal to the number of demands for all instances. Thus, demands featuring the same customer location do only occur in instances with c2c flows. To create easy and challenging instances for all targeted research areas, the number of demands is varied between 10 and 500 as the main driver of required computational effort \citep[see][]{lange2025TechRep2TCLSS}. The inner-city access \gls{tw} width is set to 2 h and 4 h, representing restricted delivery time frames \citep[see, e.g.,][]{munuzuri2012CLspain,fontaine2021scheduled}. The satellite types are chosen according to the available data and reduced amounts of road-based satellites are considered. The second-tier vehicle types are selected in line with the size of the geographic area. The demand type combinations represent three cases: inbound demands only, an inbound-outbound integration, and a desired city logistics scenario with all three types of flows. Demand sizes are generated according to the default distribution. The instances cover two extreme cases and a mixed case for the demand \gls{tw} widths, where the latter features equally many demand \glspl{tw} with short, medium, and long widths. The same holds for the customer maximum waiting time distributions. Instances without (with only long) demand \glspl{tw} and without (with only long) customer maximum waiting times are not generated. Such instances can, however, be obtained by simply ignoring the corresponding input data.
Whenever sampling is required, one replicate is created. 

Overall, the CLLIB v1 consists of 3,564 instances suitable for a variety of vehicle routing, location-routing, and service network design applications. As we show in the computational study, the number of relevant instances for a specific problem setting and solution approach are typically smaller. To support usability, we provide sub-libraries of the CLLIB with instances filtered by, e.g., the number of demands and demand types, see details in \ifbool{TR}{\ref{app:SubLibr}}{Appendix \ref{app:SubLibr}}.

\subsubsection{CLLIB-Analytics for Statistical Evaluations}

To enable a proper evaluation of effects of single instance characteristics, we provide the library extension CLLIB-Analytics, which is focused on sampling and statistical validity.
For 25, 50, and 100 demands and remaining instance generation parameter values as given in Table \ref{tab:InstParamValues}, instances are generated using a number of instance replicates parameter equal to three.
This provides the opportunity to compare otherwise identical instances with the exact same location and time characteristics, differing in only one instance characteristic value.
In total, the CLLIB-Analytics contains 93,150 instances: 32,400 each for Karlsruhe and Munich, and 28,350 for Ingolstadt.
These instances are generated independent of the CLLIB instances and, thus, qualify for out-of-sample testing of models and methods developed and tuned by using the CLLIB instances.

\subsection{Tightness of Vehicle Capacities}
\label{subsec:TightVehCap}

As mentioned in Section \ref{subsec:commonData}, for tractable instances sizes (in terms of numbers of demands), the vehicle capacities are typically no limiting factor when considering real-world-inspired package and vehicle sizes.
To thoroughly evaluate developed models and methods, however, the case of tight vehicle capacities should not be neglected. Therefore, we propose a scaling scheme for the demand and capacity data applicable to the CLLIB. 

When artificially increasing demand volumes, instance feasibility needs to be assured. This means that one package of the largest size category XL still needs to fit into one vehicle of the smallest category, here a cargo bike, and a handover at any satellite must be possible within one time period. Therefore, we multiply all demand volumes by six and all maximum demand volume per time parameters of satellites by three.  

Note that an alternative option of increasing demand volumes in a more random fashion is to use the instance generator and setting the number of demands larger than the number of customer locations. This results in instances with several demands featuring the same origin or destination location that can simply be aggregated in a post-processing step. Then, satellite capacities need to be adapted according to the resulting maximum demand volume.

%% file: CLLIB_instance_parameter_table.tex
\begin{tabular}{lw{c}{0.2\textwidth}w{c}{0.2\textwidth}w{c}{0.2\textwidth}}
	\toprule
	\textbf{City} & \textbf{Karlsruhe} & \textbf{Ingolstadt} & \textbf{Munich}  \\
	\midrule
	\textit{Number of Customer Loc./Demands} & \multicolumn{3}{c}{10, 25, 50, 100, 200, 500} \\\cmidrule{2-4}
	\textit{Planning Time Horizon} & \multicolumn{3}{c}{2 h and 4 h} \\\cmidrule{2-4}
	\textit{Satellite Type} & road, rail & road & road, rail \\\cmidrule{2-4}
	\textit{Amount of Satellites} & \multicolumn{3}{c}{as given, halved, minimum} \\\cmidrule{2-4}
	\textit{Second-tier Vehicle Type} & bike & car & car \\\cmidrule{2-4}
	\textit{Demand Types} & \multicolumn{3}{c}{e2c only, e2c (0.8) + c2e (0.2), e2c (0.65) + c2e (0.2) + c2c (0.15)} \\\cmidrule{2-4}
	\textit{Demand Size Data Source} & \multicolumn{3}{c}{distribution} \\\cmidrule{2-4}
	\textit{Demand \gls{tw} Widths Distr.} & \multicolumn{3}{c}{only short, only medium, uniform mix} \\\cmidrule{2-4}
	\textit{Customer Max. Waiting Time Distr.} & \multicolumn{3}{c}{only none, only short, uniform mix} \\\midrule
	\textbf{Number of Instances} & \textbf{1,296} & \textbf{972} & \textbf{1,296} \\
	\bottomrule
\end{tabular}

%% file: CLIG_6-CompStud.tex
In this section, we demonstrate the applicability of the CLLIB instances to two different problem settings, one concerned with scheduled service network design for tactical planning (Section \ref{subsec:tactical_planning}) and one concerned with pickups and deliveries in the city center (Section \ref{subsec:pickup_and_delivery}).

\subsection{A Scheduled Service Network Design Problem for Tactical Planning}
\label{subsec:tactical_planning}

\subsubsection{Problem Description}\label{tact_probl_descr}
The problem addresses tactical planning within the \gls{2tcl} framework. Starting from a predefined set of \glspl{dc}, multiple services, each associated with a specific vehicle type, operate along routes that connect a \gls{dc} to an ordered sequence of satellites. Final delivery from satellites to customer locations is approximated based on the distance. The task is to select an optimal set of services and assign demands accordingly, with the goal of minimizing total system cost comprising both service operating cost and assignment cost to satellites and \glspl{dc}. The problem considers only e2c demand and incorporates several capacity constraints at the satellite level. A detailed model formulation, together with the full notation, is provided in \ifbool{TR}{\ref{app:model_tactical}}{Appendix \ref{app:model_tactical}}. The formulation generalizes those of \cite{fontaine2021scheduled} and \cite{guckel2024MultiDay}.

\subsubsection{Instances}\label{tact_instances}
The numerical study considers instances with 50, 100, 200, and 500 demands for each city, focusing exclusively on inbound demands. The number and type of satellites, as well as the distribution of \glspl{tw}, are varied as described in Section \ref{subsec:CLLIB}. This yields 18 instances per demand level for Ingolstadt and 24 instances each for Munich and Karlsruhe (only road services are considered for Ingolstadt).

In addition, several problem-specific settings are introduced. Each generated service is replicated ten times within its start \gls{tw}. At every satellite visited, the first-tier vehicles operate for two periods, representing the handling time required there. To ensure feasibility with the discrete starting times of services, a minimum of 60 minutes between the release and due times of demands is enforced. The assignment cost of demands to satellites is defined as the distance divided by four. The assignment cost of demands to \glspl{dc} is based on the distance between the demand’s origin and the \gls{dc}. The distance is divided by 300 and multiplied by a random factor between 0.75 and 1.25 to capture variability in costs arising from joint delivery with other goods to the \gls{dc}. The number of urban vehicles that operate at the same time and start at the same  \gls{dc} is limited to two for the instances with up to 200 demands and to three for the instances with 500 demands.

\subsubsection{Numerical Setting and Results}\label{tact_experiments}
The model is implemented in Python 3.12 and solved with Gurobi 12. All experiments are conducted on an AMD Ryzen 9 5950X 16-Core Processor (3.40 GHz) with 128 GB RAM. For all runs, we impose a time limit of 60 minutes in Gurobi.

Table \ref{tab:results_tactical} summarizes the results grouped by the different cities and the number of demands.
It shows the average total cost, the average runtime in seconds, the number of instances solved to optimality, and the average percentage gap of the instances not solved to optimality.
As expected, increasing the number of demands results in higher objective values and longer computation times. While the model can still be solved to optimality for some instances with up to 200 demands, instances with 500 demands generally exhibit noticeable optimality gaps. A more detailed analysis of the results reveals, e.g., that assignment costs to satellites (the second-tier costs) increase substantially when the number of satellites is reduced. For the Munich instances with 500 customers, these costs rise by 41.1\% when considering the minimum number of four satellites.

\mytable[]
{Aggregate Results by City and Instance Size}
{tab:results_tactical}
{tbp}
{
\scalebox{0.95}{
    \begin{tabular}{lcrrrr}
        \toprule
        \textbf{City} & \textbf{No. Demands} & \textbf{Total Cost} & \textbf{Runtime [s]} & \textbf{Optimal} & \textbf{Gap [\%]} \\
        \midrule
        \multirow{4}{*}{Ingolstadt} & 50  & 134.42 & 0.90    & 18/18 & 0.00    \\
                                    & 100 & 252.02 & 77.83   & 18/18 & 0.00    \\
                                    & 200 & 455.79 & 3092.84 & 9/18  & 0.59 \\
                                    & 500 & 1106.05 & 3600.00 & 0/18  & 10.23 \\
        \midrule
        \multirow{4}{*}{Karlsruhe}  & 50  & 75.54  & 1.64    & 24/24 & 0.00    \\
                                    & 100 & 131.88 & 612.50  & 24/24 & 0.00    \\
                                    & 200 & 244.11 & 2271.91 & 12/24 & 0.74 \\
                                    & 500 & 592.84 & 3076.92 & 6/24  & 7.45 \\
        \midrule
        \multirow{4}{*}{Munich}     & 50  & 264.06 & 0.91    & 24/24 & 0.00    \\
                                    & 100 & 503.90 & 43.97   & 24/24 & 0.00    \\
                                    & 200 & 852.91 & 1022.90 & 20/24 & 0.20 \\
                                    & 500 & 2087.01 & 3600.00 & 0/24  & 4.59 \\
        \bottomrule
    \end{tabular}
		}
}

\subsection{A Simultaneous Pickup and Delivery Problem on the Second Tier}
\label{subsec:pickup_and_delivery}

\subsubsection{Problem Description}

We consider the \gls{sdp} as described by \cite{Bianchessi2024ResWindRed}. A set of customers $i \in N$ with demands $d_i$ for delivery and $u_i$ for pickup is given. Each customer features a service duration $s_i$ and a \gls{tw} $[e_i, l_i]$ which describes the possible start times of the service. An unlimited fleet of homogeneous capacitated vehicles operates cyclic tours from a single depot~$o$. The travel times between all pairs of locations in $N \cup \{o\}$ are known. The task is to find a set of tours with minimum total travel time so that \glspl{tw} and capacity restrictions are respected. 

\subsubsection{Instances} Since the \gls{sdp} constitutes a single-tier routing problem focused on customer locations and restricted to deliveries and pickups, only the second-tier data of the CLLIB instances with e2c and c2e demands is used. In these instances, each customer features either a delivery or pickup amount, which is defined by demand type and volume from the data. Since each customer features exactly one demand, the demand handover time defines the service time $s_i$ and the demand \glspl{tw} can be used. Travel times are immediately applicable.

Our study is based on the sub-library CLLIB-Medium with 50 and 100 demands and the urban regions Karlsruhe, Ingolstadt, and Munich. We only use instances with road-based satellites, inner-city access \glspl{tw} of 2h and 4h, as well as short, medium, and mixed demand \gls{tw} widths. Customer waiting time restrictions are not relevant. The experiments are conducted on the instances with original and scaled demand sizes (see Section \ref{subsec:TightVehCap}). Recall that the planning problem considers a single depot only, while the instance data provides several satellite locations. We, therefore, use the instances with half the number of given satellites and solve one \gls{sdp} for each satellite as depot. The instances contain 9 satellites for Munich and 5 satellites for Karlsruhe and Ingolstadt. Overall, based on 36 instances of the CLLIB-Medium, this results in 456 \gls{sdp} runs due to demand size scaling and satellite selection.

\subsubsection{Numerical Setting and Results}

We use the \gls{bpc} algorithm of \cite{Bianchessi2024ResWindRed} for solving the \gls{sdp}. The \gls{bpc} is implemented in C++ with CPLEX 22.1 as \gls{lp} solver. The experiments are run single-threaded on a high performance computing cluster using an Intel Xeon Gold 6126 with 2.6 GHz and 8GB RAM. The time limit is set to 60 minutes.
\ifbool{TR}
{The computational study examines in a first part the performance of the solution method in dependence on different instance parameter values, while we discuss the selection of a satellite as depot from economic and sustainability view points in a second part.}
{}

We first examine the performance of the \gls{bpc} algorithm depending on different instance parameter values.
Table \ref{tab:PDAvCompTimes} provides the averaged computation times grouped by city, number of demands, demand size scaling factor, and demand \gls{tw} widths distribution. As expected, instances with a larger number of demands and \gls{tw} widths are harder to solve. Considering the cities, smaller regions with shorter travel times and more symmetry in routing options are computationally more challenging.
\ifbool{TR}
{Results for different inner-city access \gls{tw} widths are not shown, since the effect of a longer access \gls{tw} width complicating the planning process is not as significant as the others.}
{}
Interestingly, more restricted instances with scaled demand sizes feature higher computation times in 14 of 18 settings. This is due to the \gls{lp} gaps being larger for these instances, resulting in larger search trees, see Table \ref{tab:PDAvLPGaps} in \ifbool{TR}{\ref{app:PDCompResults}}{Appendix \ref{app:PDCompResults}}.

\mytable[]
{Average Computation Times (in seconds) Depending on Instance Parameter Values}
{tab:PDAvCompTimes}
{tbp}
{
\scalebox{0.95}{
	\begin{tabular}{cccrrrcrrr}
		\toprule
		& && \multicolumn{3}{c}{\textbf{Regular Demand Size}} && \multicolumn{3}{c}{\textbf{Scaled Demand Size}} \\ 
		\cmidrule{4-6}\cmidrule{8-10}
		& &&\multicolumn{3}{c}{\textbf{Demand \gls{tw} Width Distr.}} && \multicolumn{3}{c}{\textbf{Demand \gls{tw} Width Distr.}} \\
		\textbf{City} & \textbf{No. of Demands} && short & medium & mixed && short & medium & mixed \\
		\midrule
		\multirow{2}{*}{Karlsruhe} & 50 && 1.68 & 878.34 & 451.27 && 2.41 & 194.31 & 100.69 \\
		& 100 && 668.08 & 2628.90 & 3600.00 && 3210.22 & 3249.50 & 3600.00 \\
		\midrule
		\multirow{2}{*}{Ingolstadt} & 50 && 0.70 & 4.97 & 2.91 && 4.64 & 15.14 & 52.54 \\
		& 100 && 26.76 & 767.26 & 1682.49 && 75.10 & 1198.29 & 3170.52 \\
		\midrule
		\multirow{2}{*}{Munich} & 50 && 0.29 & 1.52 & 17.17 && 0.52 & 2.65 & 13.13 \\
		& 100 && 10.61 & 265.96 & 402.08 && 27.48 & 895.93 & 2303.24 \\
		\bottomrule
	\end{tabular}
	}
}

We now discuss insights regarding the selection of satellites as depot from economic and sustainability view points.
Figure \ref{fig:PDAvTotalDist} shows the averaged total travel time of all vehicles resulting from different satellites selected as depot for Munich. It can be observed that the total travel time increases with a higher number of demands and less flexible \glspl{tw}. A longer inner-city access \gls{tw} causes shorter total travel times for short and mixed demand \gls{tw} widths, while this consolidation effect is not produced for medium demand \gls{tw} widths. Placing the depot at K.I. Teppichreinigung Neuhausen realizes one of the shortest total travel times for both numbers of demands and all \gls{tw} widths. While Einsteinstra\ss{}e 71 and Parkplatz Schreinerstra\ss{}e constitute reasonable choices for the 50 demand case, OBI Markt München-Westend is the second best option for the 100 demand case. All these satellites feature a central geographic location in Munich. This seems to be the most important characteristic when minimizing total travel time, since the best depot locations are rather stable against changes in the number of demands and \gls{tw} widths.

\myfigure[]
{Average Total Travel Time (in minutes) Depending on Depot Location, Number of Demands, Demand \gls{tw} Widths Distribution, and Inner-City Access \gls{tw} for Munich}
{fig:PDAvTotalDist}
{tbp}
{
	\begin{tikzpicture}
		\node[anchor=base west] at (0,0) {\includegraphics[height=50mm]{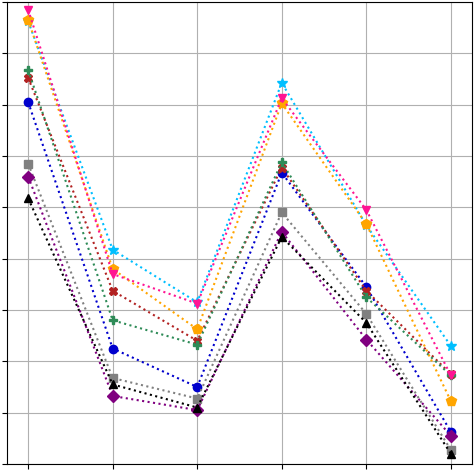}};
		\node[anchor=base west] at (5.75,0) {\includegraphics[height=50mm]{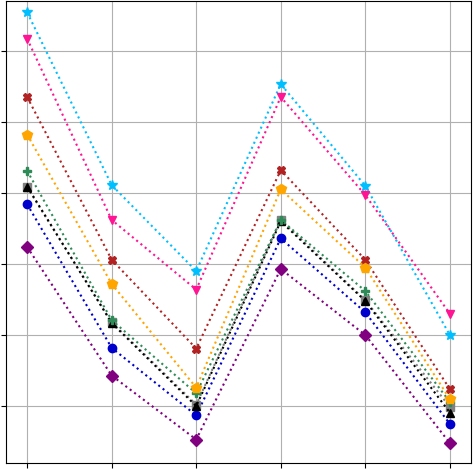}};
		\node[anchor=base west] at (11,1.5) {\includegraphics[height=27mm]{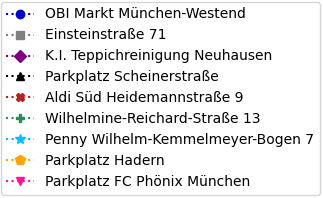}};
		% y-axes title
		\node[rotate=90] at (-0.55,2.5) {\parbox{45mm}{\centering\footnotesize Averaged Total Travel Time}};
		% y-axes ticks labels
		\node[fill=white, inner xsep=1pt,anchor=base east] at (0.1,4.91) {\tiny 620};
		\node[fill=white, inner xsep=1pt,anchor=base east] at (0.1,4.36) {\tiny 600};
		\node[fill=white, inner xsep=1pt,anchor=base east] at (0.1,3.81) {\tiny 580};
		\node[fill=white, inner xsep=1pt,anchor=base east] at (0.1,3.25) {\tiny 560};
		\node[fill=white, inner xsep=1pt,anchor=base east] at (0.1,2.71) {\tiny 540};
		\node[fill=white, inner xsep=1pt,anchor=base east] at (0.1,2.16) {\tiny 520};
		\node[fill=white, inner xsep=1pt,anchor=base east] at (0.1,1.63) {\tiny 500};
		\node[fill=white, inner xsep=1pt,anchor=base east] at (0.1,1.09) {\tiny 480};
		\node[fill=white, inner xsep=1pt,anchor=base east] at (0.1,0.53) {\tiny 460};
		\node[fill=white, inner xsep=1pt,anchor=base east] at (0.1,-0.04) {\tiny 440};	
		\node[fill=white, inner xsep=1pt,anchor=base east] at (5.85,4.38) {\tiny 1050};
		\node[fill=white, inner xsep=1pt,anchor=base east] at (5.85,3.63) {\tiny 1000};
		\node[fill=white, inner xsep=1pt,anchor=base east] at (5.85,2.85) {\tiny 950};
		\node[fill=white, inner xsep=1pt,anchor=base east] at (5.85,2.09) {\tiny 900};
		\node[fill=white, inner xsep=1pt,anchor=base east] at (5.85,1.35) {\tiny 850};
		\node[fill=white, inner xsep=1pt,anchor=base east] at (5.85,0.55) {\tiny 800};	
		% x-axes group labels
		\node at (5.5, 5.25) {\footnotesize Grouped by Number of Demands};
		\node[fill=white] at (2.7,4.65) {\scriptsize 50 demands};
		\node[fill=white] at (8.5,4.65) {\scriptsize 100 demands};	
		\node[anchor=mid east] at (0.4, -0.5) {\parbox{10mm}{\centering\footnotesize TW\\[-5pt] Widths}};
		\foreach \x/\name in {
			0.25/short,
			1.15/med,
			2.03/mix,
			2.95/short,
			3.83/med,
			4.73/mix,
			6/short,
			6.88/med,
			7.8/mix,
			8.68/short,
			9.6/med,
			10.5/mix}
		{	\draw[white] (\x,-0.8) rectangle (\x+0.375,-0.1);
			\node[rotate=90] at (\x+0.1875,-0.45) {\scriptsize \name};
		}		
		\node[anchor=mid east] at (0.4, -1.25) {\parbox{10mm}{\centering\footnotesize Access\\[-5pt] TW}};
		\draw[|-|] (0.3,-1.25) -- (2.4,-1.25);
		\node[fill=white] at (1.35,-1.25) {\footnotesize 2h};
		\draw[|-|] (3,-1.25) -- (5.1,-1.25);
		\node[fill=white] at (4.05,-1.25) {\footnotesize 4h};
		\draw[|-|] (6.05,-1.25) -- (8.15,-1.25);
		\node[fill=white] at (7.1,-1.25) {\footnotesize 2h};
		\draw[|-|] (8.73,-1.25) -- (10.88,-1.25);
		\node[fill=white] at (9.78,-1.25) {\footnotesize 4h};
	\end{tikzpicture}
}

%% file: CLIG_7-Concl.tex
% Conclusion
This paper introduces the easy-to-use and highly parametrizable open-source framework CLInstGen for generating diverse instances for \gls{2tcl} systems based on real urban data. Thereby, we summarize necessary and optional problem characteristics and highlight complicating dependencies in the instance generation. Further, we derive CLLIB, the first publicly available benchmark library for urban freight transportation with up to two tiers. Together, the framework and the library provide a standardized and extensible basis for reproducible and comparable research on different problem settings, enabling both methodological benchmarking and practical insights.

A numerical study on two problem settings from network design and vehicle routing demonstrates the versatility of the instances and produces first computational and managerial insights.

Future research can build on this work in several directions. The open-source generator provides a foundation for continuous, community-driven improvements, enabling CLLIB to evolve alongside emerging research challenges in city logistics. Based on this, the collection of real-world data can be expanded to additional cities, increasing both the diversity of generated instances and their realism. Moreover, the instance generation framework can be extended with additional transport modes and new parameters. Overall, this work lays a basis for more efficient, standardized, and generalizable city logistics research, ultimately supporting the development of more sustainable and livable cities.

%% file: CLIG_App_maps.tex
% MAPS

Figures \ref{fig:MapKA}, \ref{fig:MapING} and \ref{fig:MapMUN} show maps of the cities Karlsruhe, Ingolstadt and Munich, for which the data is provided with CLInstGen and which are used in the generation of the CLLIB. \glspl{dc} as well as road- and rail-based satellites are indicated for the urban regions, if applicable. 

\myfigure[]
{Locations of DCs and Satellites in Karlsruhe}
{fig:MapKA}
{!ht}
{
	\includegraphics[width=0.7\linewidth]{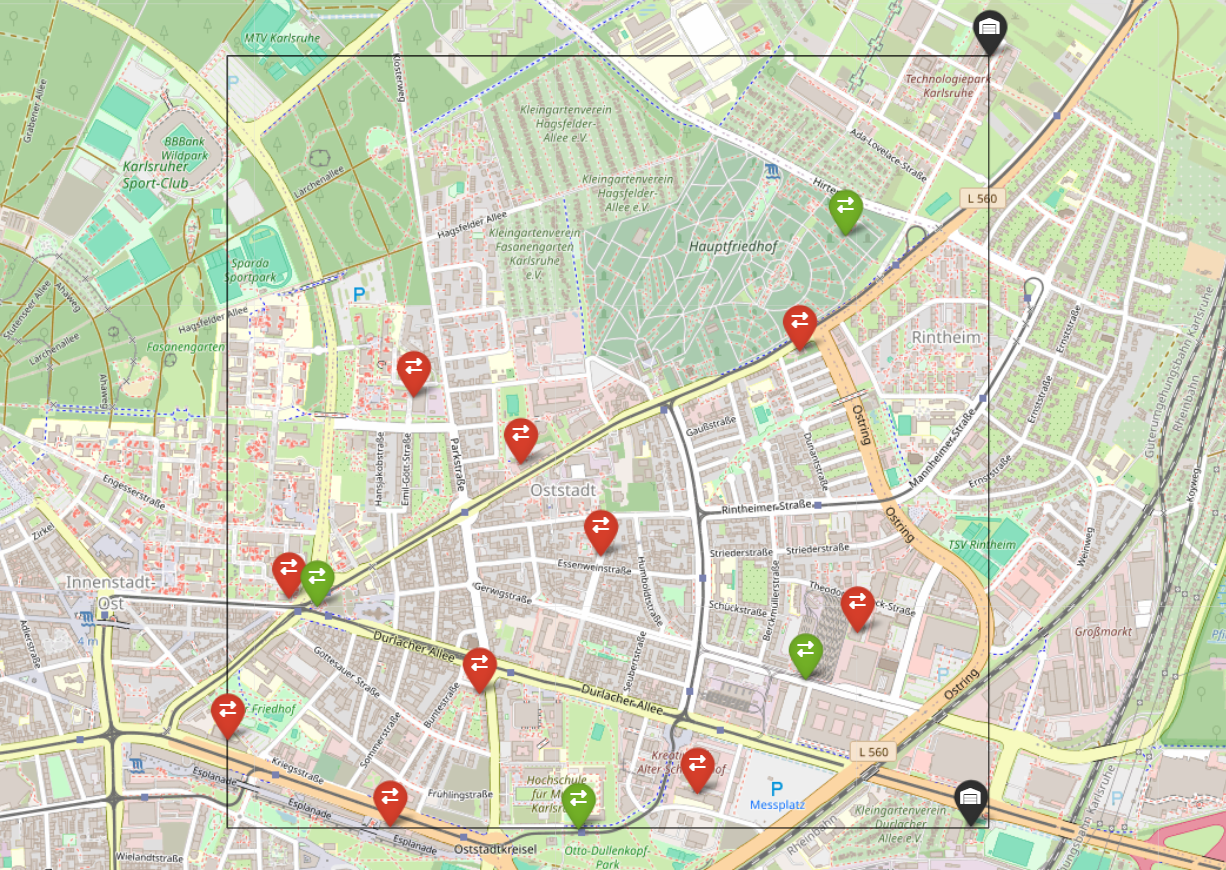}\hspace{3mm}\includegraphics[width=0.15\linewidth]{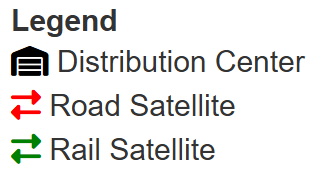}
}

\myfigure[]
{Locations of DCs and Satellites in Ingolstadt}
{fig:MapING}
{!ht}
{
	\ifbool{TR}
	{\includegraphics[width=0.83\linewidth]{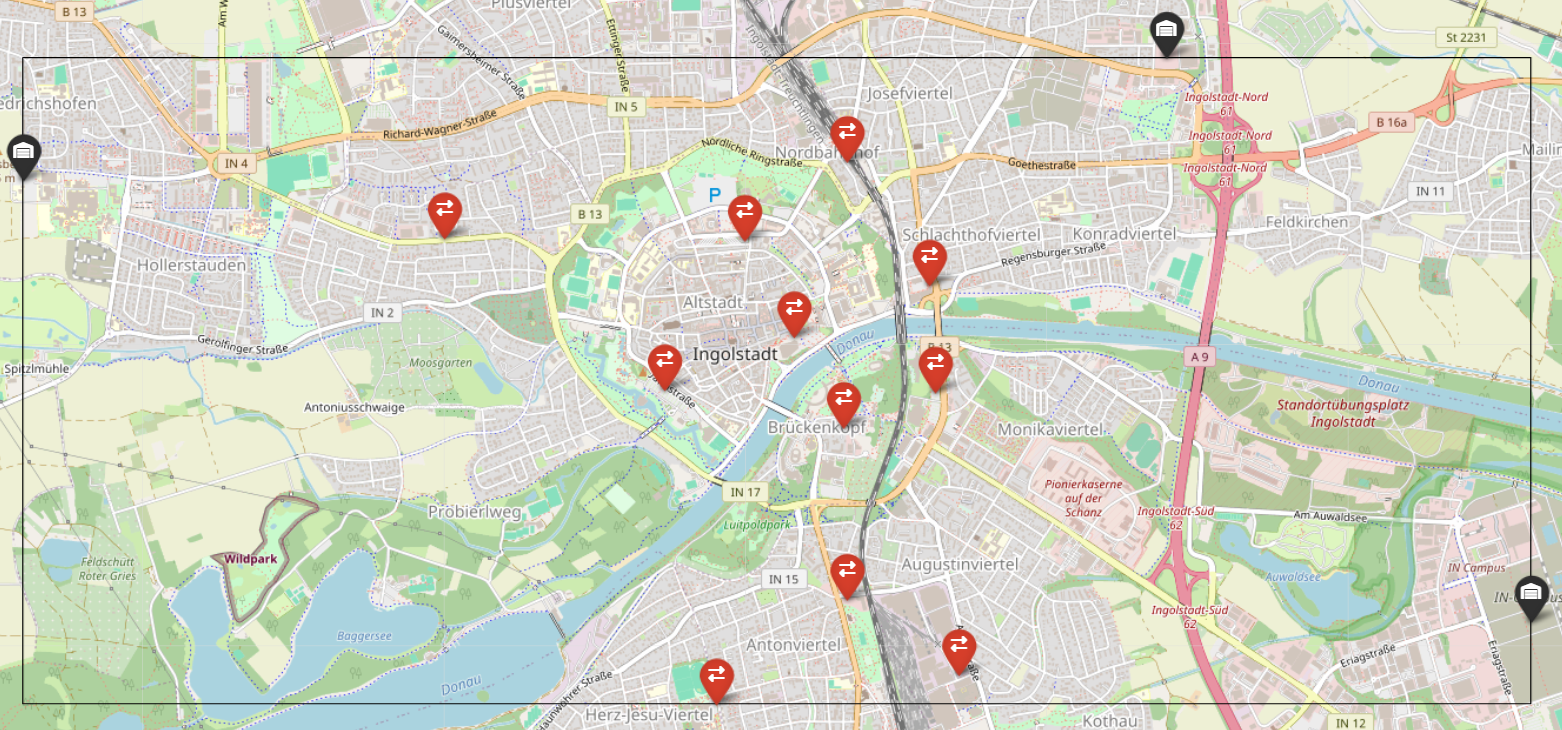}\hspace{3mm}\includegraphics[width=0.15\linewidth]{maps_legend.png}}
	{\includegraphics[width=0.84\linewidth]{map_ing.png}\hspace{3mm}\includegraphics[width=0.15\linewidth]{maps_legend.png}}
}

\myfigure[]
{Locations of DCs and Satellites in Munich}
{fig:MapMUN}
{!ht}
{
	\includegraphics[width=0.7\linewidth]{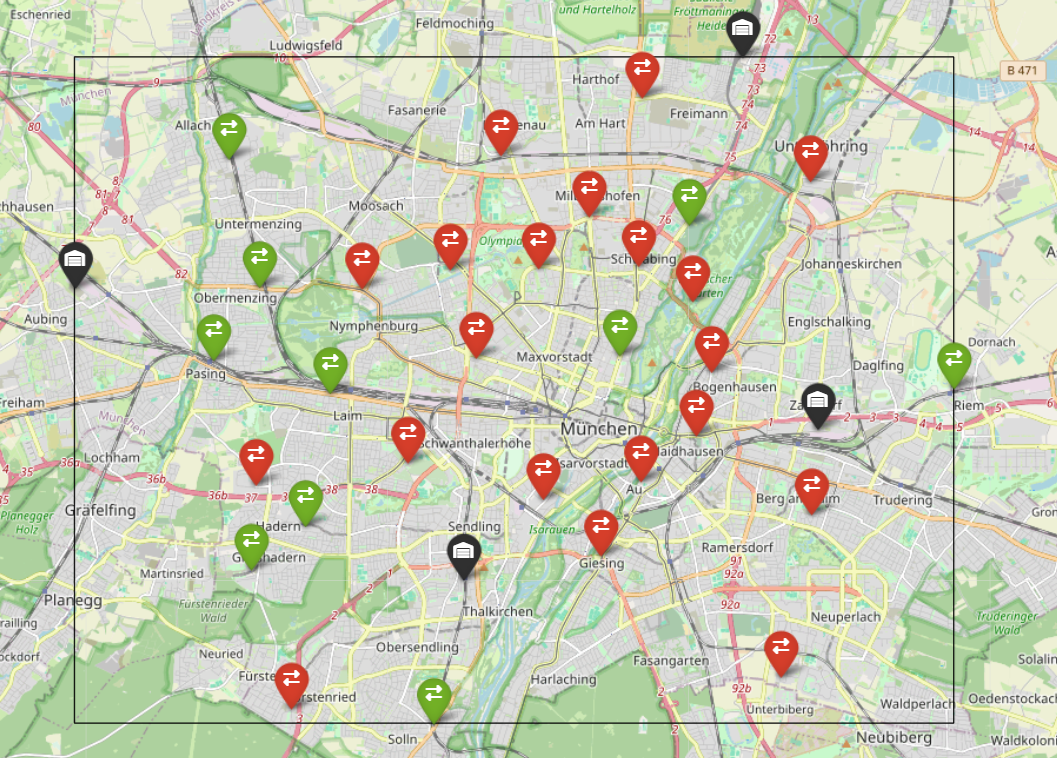}\hspace{3mm}\includegraphics[width=0.15\linewidth]{maps_legend.png}
}

%% file: CLIG_App_sublibraries.tex
% APPENDIX: CLLIB sub-libraries

To facilitate the usage of the CLLIB, sub-libraries are provided involving different subsets of instances \citep[see][]{Lange2025CLLIB}. Table \ref{tab:SubLib} summarizes the sub-libraries with the corresponding number of instances and the applied parameter value filters. In case more parameter values need to be specifically selected for an application, the CLLIB comes with a data base containing information on all instances including an instance identifier and its parameter values. Users can easily filter the required values of the parameters and obtain a list of relevant instances.

\mytable[]
{Summary of Sub-Libraries of the CLLIB}
{tab:SubLib}
{htbp}
{
	\begin{tabular}{lw{l}{0.5\textwidth}r}
		\toprule
		\textbf{Name} & \textbf{Parameter Value Filter} & \textbf{No. of Inst.} \\
		\midrule
		CLLIB-KA & city: Karlsruhe & 1,296 \\
		CLLIB-ING & city: Ingolstadt & 972 \\
		CLLIB-MUN & city: Munich & 1,296 \\
		\midrule
		CLLIB-Small & number of demands: 10, 25, 50 & 1,782 \\
		CLLIB-Medium & number of demands: 50, 100 & 1,188 \\
		CLLIB-Large & number of demands: 100, 200, 500 & 1,782 \\
		\midrule
		CLLIB-OnlyIn & demand types: e2c only & 1,188 \\
		CLLIB-InOut & demand types: e2c (0.8) + c2e (0.2) & 1,188 \\
		CLLIB-AllDemTypes & demand types: e2c (0.65) + c2e (0.2) + c2c (0.15) & 1,188 \\
		\bottomrule
	\end{tabular}
}

%% file: CLIG_App_model.tex
% DETAILS INSTANCES
\subsection{Scheduled Service Network Design Problem for Tactical Planning}
\label{app:model_tactical}

The following presents the scheduled service network design formulation of the tactical planning problem considered in Section~\ref{subsec:tactical_planning}, which is closely related with the formulations proposed in \cite{fontaine2021scheduled} and \cite{guckel2024MultiDay}. The model represents the first tier as a scheduled service network, while incorporating second tier routing costs through an approximation of the cost of assigning customer demands to satellites.

The notation is introduced in Table \ref{tab::notation}:

\mytable[]
{Notation Used in the SSND Model}
{tab::notation}
{htbp}
{
\begin{tabular}{ll}
\toprule
\textbf{Notation} & \textbf{Description} \\
\midrule
\multicolumn{2}{l}{\textbf{Sets}} \\
$\mathcal{P}$ & Set of time periods \\
$\mathcal{D}$ & Set of demands \\
$\mathcal{Z}$ & Set of satellites \\
$\mathcal{E}$ & Set of  \glspl{dc}\\
$\mathcal{R}$ & Set of urban vehicle services \\\addlinespace[0.5ex]
\multicolumn{2}{l}{\textbf{Subsets}} \\
$\mathcal{R}_{pz}$ & Set of services that operate during time period $p \in \mathcal{P}$ at satellite $z \in \mathcal{Z}$ \\
$\mathcal{R}_{pe}$ & Set of services that operate during time period $p \in \mathcal{P}$ and start at \gls{dc} $e \in \mathcal{E}$ \\
$\mathcal{R}_{dz}$ & Set of services that fulfill the time window of demand $d \in \mathcal{D}$ at satellite $z \in \mathcal{Z}$ \\\addlinespace[0.5ex]
\multicolumn{2}{l}{\textbf{Parameters}} \\
$v_{d}$ & Volume of demand $d \in \mathcal{D}$ \\
$k_{r}$ & Fixed cost of service $r \in \mathcal{R}$ \\
$f_{de}$ & Cost for assigning demand $d \in \mathcal{D}$ to \gls{dc} $e \in \mathcal{E}$ \\
$e_{r}$ & \gls{dc} from which service $r \in \mathcal{R}$ departs \\
$b_{zp}$ & Maximum demand volume that can be handled at satellite $z \in \mathcal{Z}$ during time period $p \in \mathcal{P}$ \\
$a_{zp}$ & Maximum number of urban vehicles that can operate at satellite $z \in \mathcal{Z}$ during time period $p \in \mathcal{P}$ \\
$s_{dz}$ & Approximated cost of assigning demand $d \in \mathcal{D}$ to satellite $z \in \mathcal{Z}$ \\
$u_{r}$ & Capacity of service $r \in \mathcal{R}$ \\
$n_{e}$ & Number of urban vehicles that are available at \gls{dc} $e \in \mathcal{E}$ \\\addlinespace[0.5ex]
\multicolumn{2}{l}{\textbf{Decision Variables}} \\
$y_{r}$ & Binary variable: 1 if service $r \in \mathcal{R}$ is selected, 0 otherwise \\
$x_{dzr}$ & Binary variable: 1 if demand $d \in \mathcal{D}$ is assigned to satellite $z \in \mathcal{Z}$ and service $r \in \mathcal{R}$, 0 otherwise \\
\bottomrule
\end{tabular}
}

With this notation, the mathematical formulation is as follows:

\begin{align}
\min \quad & \sum_{r \in \mathcal{R}} k_r y_r + \sum_{d \in \mathcal{D}} \sum_{z \in \mathcal{Z}} \sum_{r \in \mathcal{R}} s_{dz} x_{dzr} \label{obj} \\
\text{s.t.} \quad 
& \sum_{z \in \mathcal{Z}} \sum_{r \in \mathcal{R}_{dz}} x_{dzr} = 1 \quad && \forall d \in \mathcal{D} \label{const:assignment}\\
& \sum_{z \in \mathcal{Z}} x_{dzr} \leq y_r &&  \quad \forall d \in \mathcal{D}, \ \forall r \in \mathcal{R} \\
&\sum_{d \in \mathcal{D}} \sum_{z \in \mathcal{Z}} v_d \cdot x_{dzr} \leq u_{r} \cdot y_r && \quad \forall r \in \mathcal{R} \\
&\sum_{r \in \mathcal{R}_{pet}} y_r \leq n_{e} \quad && \forall p \in \mathcal{P},\ \forall e \in \mathcal{E},\ \forall t \in \mathcal{T} \\ 
&\sum_{d \in \mathcal{D}} \sum_{r \in \mathcal{R}_{pz}} v_d \cdot x_{dzr} \leq b_{zp} \quad && \forall z \in \mathcal{Z},\ \forall p \in \mathcal{P} \\
&\sum_{r \in \mathcal{R}_{pz}} y_r \leq a_{zp} \quad && \forall z \in \mathcal{Z},\ \forall p \in \mathcal{P} \\
& x_{dzr} \in \{0,1\} \quad && \forall d \in \mathcal{D},\ \forall z \in \mathcal{Z},\ \forall r \in \mathcal{R} \\
& y_r \in \{0,1\} \quad && \forall r \in \mathcal{R}
\end{align}

\subsection{Simultaneous Pickup and Delivery Problem on the Second Tier}
\label{app:PDCompResults}

To allow a more detailed computational analysis and to gain managerial insights for all urban regions considered, the numerical study of the \gls{sdp} has been extended.
More precisely, for all unsolved instances, we performed additional runs of the \gls{bpc} algorithm of \cite{Bianchessi2024ResWindRed} with an increased time limit of up to one week and a depth-first node-selection strategy.
With this, we obtained tight upper bounds and valid \gls{lp} bounds for all considered instances.

Table \ref{tab:PDAvLPGaps} summarizes the average \gls{lp} gaps with respect to the best known solution aggregated by city, number of demands, demand \gls{tw} width distribution, and demand size scaling. The numbers highlight the reason for the average runtimes of instances with scaled demand size being longer than those of instances with regular demand size. 
While more restricted instances tend to be easier to solve in many routing applications, the larger \gls{lp} gaps of instances with scaled demand size result in larger search trees and, thus, higher overall runtimes of the \gls{bpc} algorithm for many of the considered \gls{sdp} instances.

\mytable[]
{Average \gls{lp} Gaps Depending on Instance Parameter Values}
{tab:PDAvLPGaps}
{htbp}
{
	\begin{tabular}{cccrrrcrrr}
		\toprule
		& && \multicolumn{3}{c}{\textbf{Regular Demand Size}} && \multicolumn{3}{c}{\textbf{Scaled Demand Size}} \\ 
		\cmidrule{4-6}\cmidrule{8-10}
		& &&\multicolumn{3}{c}{\textbf{Demand \gls{tw} Width Distr.}} && \multicolumn{3}{c}{\textbf{Demand \gls{tw} Width Distr.}} \\
		\textbf{City} & \textbf{No. of Demands} && short & medium & mixed && short & medium & mixed \\
		\midrule
		\multirow{2}{*}{Karlsruhe} & 50 && 0.61\% & 0.82\% & 0.58\% && 0.66\% & 1.27\% & 1.64\% \\
		& 100 && 0.57\% & 1.06\% & 2.10\% && 1.40\% & 1.56\% & 1.84\% \\
		\midrule
		\multirow{2}{*}{Ingolstadt} & 50 && 0.48\% & 0.55\% & 0.17\% && 1.13\% & 0.95\% & 1.05\% \\
		& 100 && 0.21\% & 0.59\% & 1.00\% && 0.47\% & 0.88\% & 1.04\% \\
		\midrule
		\multirow{2}{*}{Munich} & 50 && 0.36\% & 0.68\% & 1.16\% && 0.59\% & 1.06\% & 1.09\% \\
		& 100 && 0.43\% & 0.76\% & 0.81\% && 0.79\% & 1.29\% & 1.26\% \\
		\bottomrule
	\end{tabular}
}

\afterpage{
\begin{landscape}
	\vspace*{\fill}
	\mytable[]
	{Average Total Distances for Different Depots in Munich}
	{tab:PDAvDistMunich}
	{htbp}
	{
		\begin{tabular}{lcrrrcrrrcrrrcrrr}
			\toprule
			\multicolumn{1}{c}{\textbf{Number of Demands}} && \multicolumn{7}{c}{50} && \multicolumn{7}{c}{100} \\
			\cmidrule{3-9}\cmidrule{11-17}
			\multicolumn{1}{c}{\textbf{Inner-City Access \gls{tw} Width}} && \multicolumn{3}{c}{2h} && \multicolumn{3}{c}{4h} && \multicolumn{3}{c}{2h} && \multicolumn{3}{c}{4h} \\
			\cmidrule{3-5}\cmidrule{7-9}\cmidrule{11-13}\cmidrule{15-17}
			\multicolumn{1}{c}{\textbf{Demand \gls{tw} Widths Distribution}} && \multicolumn{1}{c}{short} & \multicolumn{1}{c}{med} & \multicolumn{1}{c}{mix} && \multicolumn{1}{c}{short} & \multicolumn{1}{c}{med} & \multicolumn{1}{c}{mix} && \multicolumn{1}{c}{short} & \multicolumn{1}{c}{med} & \multicolumn{1}{c}{mix} && \multicolumn{1}{c}{short} & \multicolumn{1}{c}{med} & \multicolumn{1}{c}{mix} \\
			\midrule
			OBI Markt München-Westend && 581.0 & 485.0 & 470.0 && 553.5 & 509.0 & 452.5 && 942.0 & 841.0 & 793.5 && 918.0 & 866.0 & 787.5 \\
			Einsteinstra\ss{}e 71 && 557.0 & 473.5 & 465.5 && 538.0 & 498.5 & 445.5 && 954.0 & 859.0 & 801.0 && 931.0 & 874.5 & 799.5 \\
			K.I. Teppichreinigung Neuhausen && 552.0 & \textbf{466.5} & \textbf{461.0} && 530.5 & \textbf{488.5} & 451.0 && \textbf{912.0} & \textbf{821.5} & \textbf{776.5} && \textbf{896.5} & \textbf{850.0} & \textbf{774.0} \\
			Parkplatz Scheinerstra\ss{}e && \textbf{543.5} & 471.0 & 462.0 && \textbf{528.5} & 495.0 & \textbf{444.0} && 954.5 & 858.5 & 800.0 && 930.0 & 874.0 & 795.5 \\
			Aldi Süd Heidemannstra\ss{}e 9 && 590.5 & 507.5 & 488.0 && 555.5 & 507.5 & 475.0 && 1017.5 & 903.0 & 840.5 && 966.0 & 903.0 & 812.0 \\
			Wilhelmine-Reichard-Stra\ss{}e 13 && 593.5 & 496.0 & 486.5 && 557.5 & 505.0 & 474.5 && 965.5 & 860.5 & 809.0 && 931.0 & 881.0 & 803.0 \\
			Penny Wilhelm-Kemmelmeyer-Bogen 7 && 612.5 & 523.5 & 503.0 && 588.5 & 533.0 & 486.0 && 1077.0 & 955.5 & 895.0 && 1026.5 & 955.0 & 850.0 \\
			Parkplatz Hadern && 613.0 & 516.0 & 492.5 && 580.5 & 533.5 & 464.5 && 991.0 & 886.0 & 813.0 && 952.5 & 897.5 & 805.0 \\
			Parkplatz FC Phönix München && 617.0 & 514.0 & 502.5 && 582.5 & 539.0 & 474.5 && 1058.5 & 931.0 & 882.0 && 1017.5 & 948.5 & 865.0 \\
			\bottomrule
		\end{tabular}
	}
\vspace*{\fill}\mbox{}
\end{landscape}
}